\documentclass[11pt]{article}
\usepackage{txfonts}
\usepackage{color,epsfig}
\usepackage{mathrsfs}
\usepackage{amssymb}
\usepackage{amsfonts,graphicx,psfrag}

\newtheorem{Th}{Theorem}[section]

\newtheorem{lemma}{Lemma}[section]
\newtheorem{theorem}[lemma]{Theorem}

\newtheorem{proposition}[lemma]{Proposition}

\newtheorem{remark}[lemma]{Remark}
\let\lutzremark=\remark
\def\remark{\lutzremark\normalfont}

\newtheorem{notation}[lemma]{Notation}
\newtheorem{example}[lemma]{Example}

\def\be{\begin{equation}}
\def\ee{\end{equation}}
\def\bea{\begin{eqnarray}}
\def\eea{\end{eqnarray}}
\def\bes{\begin{eqnarray*}}
\def\ees{\end{eqnarray*}}

\def\nn{\nonumber}
\def\<{\langle}
\def\>{\rangle}
\def\lb{\label}
\def\bs{\setminus}
\def\pt{\partial}

\def\R{{\bf R}}

\def\Z{{\bf Z}}
\def\N{{\bf N}}

\def\SO{{\rm SO}}

\def\aa{{\alpha}}
\def\bb{{\beta}}
\def\ga{{\gamma}}

\def\th{{\theta}}

\def\ep{{\epsilon}}
\def\lm{{\lambda}}

\def\sg{{\sigma}}

\def\diag{{\rm diag}}

\def\hb{\vrule height0.18cm width0.14cm $\,$}

\def\td#1{\tilde{#1}}

\title{The reduction on the linear stability of elliptic Euler-Moulton solutions of
   the $n$-body problem to those of $3$-body problems}

\author{Qinglong Zhou$^{1} $\thanks{Partially supported by NSFC (No.11501330, No.11425105) of China.
           E-mail:zhouqinglong@sdu.edu.cn}\quad
        Yiming Long$^{2} $\thanks{Partially supported by NSFC (No. 11131004), LPMC of MOE of China,
           Nankai University, and BAICIT at Capital Normal University. E-mail: longym@nankai.edu.cn}, \\ \\
$^{1}$ School of Mathematics\\Shandong University, Jinan 250100, Shandong, China\\
$^{2}$ Chern Institute of Mathematics and LPMC\\Nankai University, Tianjin 300071, China\\
}

\begin{document}

\maketitle

\begin{abstract}
{In this paper, we consider the elliptic collinear solutions of the classical $n$-body problem, where
the $n$ bodies always stay on a straight line, and each of them moves on its own elliptic orbit with
the same eccentricity. Such a motion is called an elliptic Euler-Moulton collinear solution. Here we
prove that the corresponding linearized Hamiltonian system at such an elliptic Euler-Moulton collinear
solution of $n$-bodies splits into $(n-1)$ independent linear Hamiltonian systems, the first one is the
linearized Hamiltonian system of the Kepler $2$-body problem at Kepler elliptic orbit, and each of the
other $(n-2)$ systems is the essential part of the linearized Hamiltonian system at an elliptic Euler
collinear solution of a $3$-body problem whose mass parameter is modified. Then the linear stability of
such a solution in the $n$-body problem is reduced to those of the corresponding elliptic Euler collinear
solutions of the $3$-body problems, which for example then can be further understood using numerical
results of Martin\'ez, Sam\`a and Sim\'o in \cite{MSS1} and \cite{MSS2} on $3$-body Euler solutions in
2004-2006. As an example, we carry out the detailed derivation of the linear stability for an elliptic
Euler-Moulton solution of the $4$-body problem with two small masses in the middle.}
\end{abstract}



\renewcommand{\theequation}{\thesection.\arabic{equation}}

\setcounter{equation}{0}
\section{Introduction and main results}
\label{sec:1}

When one considers a system of $n$ bodies including the Earth, the Moon and $(n-2)$ space stations in
the middle, one tries to find places for these space stations so that they can be easily put there and
easily taken away. When $n=3$, by the linear stability study it is well-known that such a middle place
should be the Euler point, because at such a point the essential part of the linearized Hamiltonian
system possesses two pairs of Floquet multipliers with suitable masses and eccentricity, one of which
is elliptic and the other is hyperbolic. This paper is devoted to study the problem for general $n\ge 3$,
and in fact here we prove that the study on such an $n$-body problem can be reduced to those of $(n-2)$
related $3$-body problems.

Recall that in the classical $3$-body problem with three positive masses, a special solution was found
by L. Euler in \cite{Euler} of 1767. In this motion, the $3$ bodies form always a collinear central
configuration at any time in a fixed plane and each body runs along a special Keplerian elliptic orbit
about the center of mass of the $3$ bodies with the same eccentricity $e\in [0,1)$. Then F. Moulton in
\cite{Mou} of 1910 proved that for every ordering of $n$ positive masses, there exists a unique
collinear central configuration of $n$-bodies. After them in general, for the classical $n$-body problem
we call a solution {\it elliptic Euler-Moulton homographic motion} of $n$-bodies ({\it EEM} for short
below), if the $n$ bodies always form a collinear central configuration and each body travels along a
specific Keplerian elliptic orbit about the center of mass of the system with the same eccentricity.
Specially when $e=0$, the $n$ bodies run circularly around the center of mass with the same angular
velocity, which are called {\it Euler-Moulton relative equilibria} traditionally.

Given $n$ positive masses $m=(m_1, \ldots, m_n)\in (\R^+)^n$ on $n$ points $q=(q_1,\ldots,q_n)\in (\R^2)^n$
respectively. According to Newton's gravitation law, their motion is governed by the system,
\be   m_i\ddot{q}_i=\frac{\partial U(q)}{\partial q_i}, \qquad {\rm for}\quad i=1, 2, \ldots, n, \lb{1.1}\ee
where $U(q) = \sum_{1\leq i<j\leq n}\frac{m_im_j}{|q_i-q_j|}$ is the potential function and $|\cdot|$
denotes the norm of vectors in $\R^2$.

Let
$$  \hat{\mathcal {X}}:=\left\{q=(q_1,q_2,\ldots,q_n)\in (\R^2)^n\,\,\left|\,\,
       \sum_{i=1}^n m_iq_i=0,\,\,q_i\neq q_j,\,\,\forall i\neq j \right. \right\}.  $$
Then critical points of the action functional
$$ \mathcal{A}(q)=\int_{0}^{2\pi}\left[\sum_{i=1}^n\frac{m_i|\dot{q}_i(t)|^2}{2}+U(q(t))\right]dt $$
defined on the space $W^{1,2}(\R/2\pi\Z,\hat{\mathcal {X}})$ correspond to $2\pi$-periodic solutions of the
system (\ref{1.1}) one-to-one. To transform (\ref{1.1}) to a Hamiltonian system, we let
$p=(p_1, \ldots, p_n)$ with $p_i=m_i\dot{q}_i\in\R^2$ for $1\le i\le n$ and obtain
\be \dot{p}_i=-\frac{\partial H}{\partial q_i},\,\,\dot{q}_i
  = \frac{\partial H}{\partial p_i},\qquad {\rm for}\quad i=1,2,\ldots, n,  \lb{1.2}\ee
where the Hamiltonian function is given by
\be H(p,q) = \sum_{i=1}^n\frac{|p_i|^2}{2m_i} - U(q).  \lb{1.3}\ee

It is well-known (cf. \cite{MSS1}, \cite{MSS2}) that the linear stability of an EEM solution
of the $3$-body problem with masses $m=(m_1,m_2,m_3)\in (\R^+)^3$ is determined by the
eccentricity $e\in [0,1)$ and the mass parameter
\be  \bb = \frac{m_1(3x^2+3x+1)+m_3x^2(x^2+3x+3)}{x^2+m_2[(x+1)^2(x^2+1)-x^2]},    \lb{1.4}\ee
where $x$ is the unique positive solution of the Euler quintic polynomial equation
\be (m_3+m_2)x^5+(3m_3+2m_2)x^4+(3m_3+m_2)x^3-(3m_1+m_2)x^2-(3m_1+2m_2)x-(m_1+m_2)=0, \lb{1.5}\ee
and the three bodies form a central configuration of $m$, which are denoted by $q_1=0$,
$q_2=(x\alpha,0)^T$ and $q_3=((1+x)\alpha,0)^T$ with $\alpha=|q_2-q_3|>0$, $x\alpha=|q_1-q_2|$.

In this paper we prove that the linear stability problem of the EEM in the $n$-body case for every
integer $n\ge 3$ can be in fact reduced to the linear stabilities of $(n-2)$ related EEM of $3$-body
cases. More precisely, based on the central configuration coordinate method of K. Meyer and D. Schmidt
in \cite{MS}, we reduce the linear stability of the $n$-body EEM to two parts symplectically, one of
which is the same as that of the Kepler solutions, and the other is a $4(n-2)$-dimensional Hamiltonian
system whose fundamental solution is the essential part for the linear stability of the EEM of
$n$-bodies. Then we prove that this essential part is the sum of $(n-2)$ independent linear Hamiltonian
systems, each of which is the essential part of the linearized Hamiltonian system of some EEM of a
related $3$-body problem.

To describe our main reduction result more precisely, given positive masses
$m=(m_1,m_2,\ldots,m_n)\in (\R^+)^n$, let $a=(a_1, \ldots, a_n)$ be the unique $n$-body collinear
central configuration of $m$ with $a_i=(a_{ix},0)^T$ for $1\le i\le n$ which satisfies $a_{ix}<a_{jx}$
if $i<j$. Without lose of generality, we normalize the masses by
\be   \sum_{i=1}^n m_i=1, \lb{1.6}\ee
and normalize the positions $a_i$ with $1\le i\le n$ by
\be  \sum_{i=1}^n m_ia_i=0, \quad{\rm and}\quad \sum_{i=1}^n m_ia_i^2=2I(a)=1. \lb{1.7}\ee
Moreover, we define
\be  \mu = U(a) = \sum_{1\le i<j\le n}\frac{m_im_j}{|a_i-a_j|},  \lb{1.8}\ee
and
\be  \tilde{M}=\diag(m_1, \ldots, m_n).    \lb{1.9}\ee
Let $B=(B_{ij})=\frac{1}{2}U''(a)$ be the Hessian of $U(q)$ at the collinear central configuration $q=a$
which is an $n\times n$ symmetric matrix given by
\bea
B_{ij} &=& \frac{m_im_j}{|a_i-a_j|^3},\qquad{\rm if}\;i\ne j,\;1\le i,j\le n, \lb{1.10}\\
B_{ij} &=& -\sum_{1\le j\le n\atop j\ne i}\frac{m_im_j}{|a_i-a_j|^3},\qquad{\rm if}\;i=j,\;1\le i\le n. \lb{1.10+1}\eea
We let
\be  D = \mu I_n + \td{M}^{-1}B.  \lb{1.11}\ee

Then the following lemma is crucial for our study, whose proof is due to C. Conley according to
F. Pacella (\cite{Pac}, 1987) and R. Moeckel (\cite{Moe1} of 1990 as well as \cite{Moe2} of 1994).
For reader's conveniences, a sketch of this proof will be given in the Appendix of this paper below
following \cite{Pac}, \cite{Moe1} and \cite{Moe2}.

\begin{lemma}\label{L1.1}
The $n\times n$ matrix $D$ possesses a simple eigenvalue $\lm_1=\mu>0$ and a second eigenvalue
$\lm_2=0$. The other $n-2$ eigenvalues of $D$ besides $\mu$ and this $0$ are non-positive.
Consequently they satisfy
\be   \lm_1 > \lm_2=0 \ge \lm_3 \ge \cdots \ge \lm_n.  \lb{1.12}\ee
\end{lemma}

Then we define
\be  \bb_i = -\frac{\lm_{i+2}}{\mu} \ge 0, \qquad \forall\; 1\le i\le n-2.  \lb{1.13}\ee
Based on these $\bb_i$s, our main result of this paper is the following

\begin{theorem}\label{T1.2}
In the planar $n$-body problem with given masses $m=(m_1,m_2,\ldots,m_n)\in (\R^+)^n$, denote the
EEM with eccentricity $e\in [0,1)$ for $m$ by $q_{m,e}(t)=(q_1(t), q_2(t), \ldots,q_n(t))$. Then
the linearized Hamiltonian system at $q_{m,e}$ is reduced into the sum of $(n-1)$ independent
Hamiltonian systems, the first one is the linearized system of the Kepler $2$-body problem at
the corresponding Kepler orbit, and the $i$-th part of the other $(n-2)$ parts with $1\le i\le n-2$
is the essential part of the linearized Hamiltonian system of some EEM of a $3$-body problem with the
original eccentricity $e$ and the mass parameter $\bb_i$ given by (\ref{1.13}) instead of that
$\bb$ given by (\ref{1.4}).
\end{theorem}

\begin{remark}\label{R1.3}
(i) J. Liouville first observed in \cite{Lio} of 1842 that the Moon stays always on the straight
line passing through the centers of the Sun and the Earth and on the opposite side of the Sun with
respect to the Earth, i.e., the Moon always enlightens the Earth during the nights, is impossible
due to the instability of such a configuration. According to R. Moeckel (cf. p.300, \cite{Moe2}) of
1994, the stability analysis of collinear relative equilibria can be attributed to M. Andoyer
\cite{And} in 1906 and M. Meyer \cite{Mey} in 1933. Subsequent studies on the linear stability of
EEMs can be found in \cite{MS} of K. Meyer and D. Schmidt in 2005, \cite{MSS}, \cite{MSS1} and
\cite{MSS2} of R. Mart\'{\i}nez, A. Sam\`{a} and C. Sim\'{o} in 2004-2006, and the recent preprints
\cite{ZL} of Q. Zhou and Y. Long, and \cite{HO} of X. Hu and Y. Ou. Researches on Lagrangian equilateral
triangle elliptic solutions (cf. \cite{Lag}) and related topics were done by M. Gascheau (\cite{Ga},
1843), E. Routh (\cite{R2}, 1875), J. Danby (\cite{Dan}, 1964), R. Moeckel (\cite{Moe3}, 1995), G.
Roberts (\cite{R1}, 2002), X. Hu and S. Sun (\cite{HS}, 2010), and X. Hu, Y. Long and S. Sun
(\cite{HLS}, 2014).

(ii) Based on our above reduction theorems, the numerical results obtained by R. Mart\'{\i}nez, A.
Sam\`{a} and C. Sim\'{o} in \cite{MSS1} and \cite{MSS2} for $3$-body Euler solutions can be applied to
get the linear stability of the $n$-body elliptic Euler-Moulton collinear solutions using our formula of
$\bb_i$s in (\ref{1.13}) for any positive integer $n\ge 3$. The theoretical linear stability results
on $3$-body EEM obtained in papers \cite{ZL} and \cite{HO} can also be applied too.

(iii) It may be worth to point out that the proof of our reduction Theorem \ref{T1.2} is based upon
the results of \cite{MS} of 2005, and is independent of the results and their proofs in papers
\cite{MSS1}, \cite{MSS2}, \cite{ZL} and \cite{HO} for the $3$-body case.
\end{remark}

In the Section 2 of this paper we focus on the proof of Theorem \ref{T1.2}. In Section 3,
we study a special example of a collinear $4$-body problem with two small masses in the
middle. The two corresponding mass parameters $\bb_1$ and $\bb_2$ in (\ref{1.13}) are calculated
explicitly there, and hence their linear stability can be determined numerically using results
in \cite{MSS1} and \cite{MSS2} of 2004-2006 for example. It is interesting to see that when the
masses of the two middle particles tend to $0$, the effect of both of them does not disappear.
In the Appendix, a sketch of the proof of Lemma 1.1 is given.

\setcounter{equation}{0}
\section{Reduction from the collinear $n$-body problem to $(n-2)$ collinear $3$-body problems}\label{sec:2}

In their paper \cite{MS} of 2005, K. Meyer and D. Schmidt introduced the central configuration coordinates
for a class of periodic solutions of the $n$-body problem. Our study on the EEM solutions of $n$-bodies
is based upon their method. Here the key point is that we found the reduction of the linear stability
of the $n$-body EEM problem to those of $(n-2)$ three body problems. This reduction needs more techniques
for the $n$ body case.

As in Section 1, for the given masses $m=(m_1,m_2,\ldots,m_n)\in (\R^+)^n$ satisfying (\ref{1.6}),
suppose the $n$ particles are all on the $x$-axis with $a_1=(a_{1x},0)^T$, $a_2=(a_{2x},0)^T$,
$\ldots$, $a_n=(a_{nx},0)^T$ satisfying $a_{ix}<a_{jx}$ if $i<j$. In this section we always denote by
$a=(a_1,\ldots,a_n)$ the unique collinear central configuration for the mass $m$ determined by
\cite{Mou}. Using normalization and notations (\ref{1.6})-(\ref{1.9}), we have
\be \sum_{j=1,j\ne i}^n\frac{m_j(a_{jx}-a_{ix})}{|a_{jx}-a_{ix}|^3}
          = \frac{U(a)}{2I(a)}a_{ix}=\mu a_{ix}.    \label{eq.of.cc}\ee
Based on the matrix $B$ of (\ref{1.10})-(\ref{1.10+1}), besides $D$ we further define
\be  \tilde{D} = \mu I_n+\tilde{M}^{-1/2}B\tilde{M}^{-1/2}
     = \tilde{M}^{1/2}D\tilde{M}^{-1/2}. \label{tilde.D}\ee
where $\mu$ is given by (\ref{1.8}).

Since $\td{D}$ is symmetric, all its eigenvalues are real, which are denoted by $\lm_1=\mu$,
$\lm_2=0$, $\lm_3$, $\ldots$, $\lm_n$ with corresponding eigenvectors $\td{v}_1=\td{M}^{1/2}v_1$,
$\td{v}_2=\td{M}^{1/2}v_2$, $\td{v}_3$, $\ldots$, $\td{v}_n$. Moreover, we can suppose that
$\td{v}_1$, $\td{v}_2$, $\ldots$, $\td{v}_n$ form an orthonormal basis of $\R^n$.

Letting $v_i=\tilde{M}^{-1/2}\tilde{v}_i$ for $3\le i\le n$, we have
$$ Dv_i=\tilde{M}^{-1/2}\tilde{D}\tilde{M}^{1/2}(\tilde{M}^{-1/2}\tilde{v}_i)=\tilde{M}^{-1/2}\tilde{D}\tilde{v}_i
    =\tilde{M}^{-1/2}\lambda_i\tilde{v}_i=\lambda_i v_i.  $$
Thus $v_i$ is the eigenvector of $D$ belonging to its eigenvalue $\lambda_i$. Moreover, by the
orthonormal basis property of $\tilde{v}_1$, $\tilde{v}_2$, $\ldots$, $\tilde{v}_n$, we have
\be
v_i^T\tilde{M}v_j=\tilde{v}_i^T\tilde{v}_j=\delta_i^j,\qquad \forall\;1\le i,j\le n.  \label{vMv}
\ee

Denote the eigenvector $v_i$ belonging to the eigenvalue $\lm_i$ of the matrix $D$ by
$v_i^T=(b_{1i},b_{2i},\ldots,b_{ni})$ for $3\le i\le n$, i.e.,
\be  D(b_{1k},b_{2k},\ldots,b_{nk})^T=\lambda_k(b_{1k},b_{2k},\ldots,b_{nk})^T, \qquad 3\le k\le n. \lb{bi}\ee

Then it yields
\be \mu b_{ik}-\sum_{j=1,j\ne i}^n\frac{m_j(b_{ik}-b_{jk})}{|a_i-a_j|^3}=\lambda_k b_{ik},\quad 1\le i\le n. \ee
Let
\be
F_{ik}=\sum_{j=1,j\ne i}^n\frac{m_im_j(b_{ik}-b_{jk})}{|a_i-a_j|^3},\quad 1\le i\le n,3\le k\le n, \label{F_ki}
\ee
then we have
\be
F_{ik}=(\mu-\lambda_k)m_ib_{ik}. \label{Fi.bi}
\ee
Moreover, we have
\be
\sum_{i=1}^nF_{ik}b_{ik}=\sum_{i=1}^n(\mu-\lambda_k)m_ib_{ik}^2
    =\mu-\lambda_k=\mu(1+\bb_{k-2}), \label{sum.of.Fi.bi}\ee
where in the last equality, we used (\ref{1.13}).

Now as in p.263 of \cite{MS}, we define
\begin{equation}\label{PQYX}
P=\left(\matrix{p_1\cr p_2\cr \ldots\cr p_n}\right),
\quad
Q=\left(\matrix{q_1\cr q_2\cr \ldots\cr q_n}\right),
\quad
Y=\left(\matrix{G\cr Z\cr W_1\cr \ldots\cr W_{n-2}}\right),
\quad
X=\left(\matrix{g\cr z\cr w_1\cr \ldots\cr w_{n-2}}\right),
\end{equation}
where $p_i$, $q_i$ with $i=1,2,\ldots,n$, $G$, $Z$, $W_i$ with $1\le i\le n-2$,
$g$, $z$, and $w_i$ with $1\le i\le n-2$ are all column vectors in $\R^2$. We make
the symplectic coordinate change
\be\lb{transform1}  P=A^{-T}Y,\quad Q=AX,  \ee
where the matrix $A$ is constructed as in the proof of Proposition 2.1 in \cite{MS}.
More precisely, the matrix $A\in {\bf GL}(\R^{2n})$ is given by
\begin{equation}
A=
\left(
\matrix{
I\quad A_1\quad B_{13}\quad \ldots\quad B_{1n}\cr
I\quad A_2\quad B_{23}\quad \ldots\quad B_{2n}\cr
\ldots\quad \ldots\quad \ldots\quad \ldots\quad \ldots\cr
I\quad A_n\quad B_{n3}\quad \ldots\quad B_{nn}
}
\right),
\end{equation}
where each $A_i$ is a $2\times2$ matrix given by
\be A_i = (a_i, Ja_i)=\left(\matrix{a_{ix}\quad 0\cr 0\quad a_{ix}}\right)=a_{ix}I_2. \label{Aa}\ee
Let
\be B_{ki} = \left(\matrix{b_{ki}\quad 0\cr 0\quad b_{ki}}\right)=b_{ki}I_2.  \label{Bb}\ee
Then $A^TMA=I_{2n}$ holds (cf. (13) in p.263 of \cite{MS}).

As in Theorem 2.1 on pp.261-262, setting $G=g=0$ to fix the center of mass at the origin as in p.271
of \cite{MS}, after the transform (\ref{transform1}) the Hamiltonian function of the $n$-body problem
in the new variables becomes
\be H(Z,W_1,\ldots,W_{n-2},z,w_1,\ldots,w_{n-2}) = K(Z, W_1,\ldots,W_{n-2}) - U(z,w_1,\ldots,w_{n-2})   \lb{HaFun}\ee
where the kinetic energy satisfies
\be K = \frac{1}{2}(|Z|^2+|W_1|^2+\ldots+|W_{n-2}|^2), \ee
and the potential function satisfies
\be  U(z,w_1,\ldots,w_{n-2}) = \sum_{1\le i<j\le n}U_{ij}(z,w_1,\ldots,w_{n-2}), \label{U}\ee
with
\bea
U_{ij}(z,w_1,\ldots,w_{n-2}) &=& \frac{m_im_j}{d_{ij}(z,w_1,\ldots,w_{n-2})}, \label{U_ij}\\
d_{ij}(z,w_1,\ldots,w_{n-2}) &=& |(A_i-A_j)z+\sum_{k=3}^{n}(B_{ik}-B_{jk})w_{k-2}|  \nn\\
 &=& |(a_{ix}-a_{jx})z+\sum_{k=3}^{n}(b_{ik}-b_{jk})w_{k-2}|,  \label{d_{ij}}\eea
where we have used (\ref{Aa}) and (\ref{Bb}). Recall that each $Z$, $W_i$, $z$, $w_i$ with $1\le i\le n-2$ is
a vector in $\R^2$. Here $z=z(t)$ is the Kepler elliptic orbit given through the true anomaly $\th=\th(t)$,
\be  r(\th(t)) = |z(t)| = \frac{p}{1+e\cos\th(t)},  \lb{rTh}\ee
where $p=a(1-e^2)$ and $a>0$ is the latus rectum of the ellipse (\ref{rTh}).

As in pp.271-273 of \cite{MS}, we have the following proposition.

\begin{proposition}\label{P2.1}
There exists a symplectic coordinate change
\be   \xi = (Z,W_1,\ldots,W_{n-2},z,w_1,\ldots,w_{n-2})^T
    \;\mapsto\; \bar{\xi}
    = ( \bar{Z}, \bar{W}_1,\ldots, \bar{W}_{n-2}, \bar{z}, \bar{w}_1,\ldots, \bar{w}_{n-2})^T, \lb{P21-1}\ee
such that using the true anomaly $\th$ as the variable the resulting Hamiltonian function of the
$n$-body problem is given by
\bea
&& H(\theta,\bar{Z},\bar{W_1},\ldots,\bar{W}_{n-2},\bar{z},\bar{w_1},\ldots,\bar{w}_{n-2})  \nn\\
&&\qquad = \frac{1}{2}(|\bar{Z}|^2 + \sum_{k=1}^{n-2}|\bar{W_k}|^2)
      + (\bar{z}\cdot J\bar{Z} + \sum_{k=1}^{n-2}\bar{w_k}\cdot J\bar{W_k}) \nn\\
&&\qquad\quad + \frac{p-r(\theta)}{2p}(|\bar{z}|^2+\sum_{k=1}^{n-2}|\bar{w_k}|^2)
          - \frac{r(\theta)}{\sigma}U(\bar{z},\bar{w_1},\ldots,\bar{w}_{n-2}), \label{P21-2}\eea
where $J = \left(\matrix{0 &-1\cr
                         1 & 0\cr}\right)$, $r(\th)=\frac{p}{1+e\cos\th}$,
$\mu$ is given by (\ref{1.8}), $\sg=(\mu p)^{-1/4}$ and $p$ is given in (\ref{rTh}).
\end{proposition}

\begin{remark}\lb{correction} Proposition \ref{P2.1} is a modified version of Lemma 3.1 of \cite{MS}
in our case of $n$-bodies. As pointed out in Section 11 of \cite{Lon}, in the $3$-body case, the
$\sg$ in (\ref{P-02}) given by $\sg=p\beta^3$ in the original computation on line 9 of p.273 in
\cite{MS} is incorrect, and should be corrected to $\sg=(\mu p)^{-1/4}$. Note also that in the line
11 of p.273 in \cite{MS}, the stationary solution $(0,1,0,0,1,0,0,0)^T$ is not correct too and should
be corrected to $(0,\sg,0,0,\sg,0,0,0)^T$ as in \cite{Lon}, and in general it may not be
possible to have $\sg=1$.
\end{remark}

{\bf Proof of Proposition \ref{P2.1}.} Because of reasons mentioned in this remark, for reader's
conveniences, we give the complete details of the proof of this proposition below.

Following the proof of Lemma 3.1 of \cite{MS}, we carry the coordinate changes in four steps.

{\bf Step 1.} {\it Rotating coordinates via the matrix $R(\th(t))$ in time $t$.}

We change first the coordinates $\xi$ to
\be \hat{\xi}=(\hat{Z}, \hat{W}_1, \ldots, \hat{W}_{n-2}, \hat{z}, \ldots, \hat{w}_1, \hat{w}_{n-2})^T
                  \in (\R^2)^{n-1}, \lb{P-01}\ee
which rotates with the speed of the true anomaly. The transformation matrix is given by the rotation
matrix $R(\th) = \left(\matrix{\cos\th &-\sin\th\cr
                               \sin\th & \cos\th\cr}\right)$. The generating function of this
transformation is given by
\be  \hat{F}(t, Z, W_1, \ldots, W_{n-2}, \hat{z}, \hat{w}_1, \ldots, \hat{w}_{n-2})
       = -Z\cdot R(\th)\hat{z} - \sum_{i=1}^{n-2}W_i\cdot R(\th)\hat{w}_i,  \lb{P-02}\ee
and for $1\le i\le n-2$ the transformation is given by
\bea
z = -\frac{\pt \hat{F}}{\pt Z} = R(\th)\hat{z}, && \hat{Z} = -\frac{\pt \hat{F}}{\pt \hat{z}} = R(\th)^TZ,  \lb{P-03}\\
w_i = -\frac{\pt \hat{F}}{\pt W_i} = R(\th)\hat{w}_i, && \hat{W}_i = -\frac{\pt \hat{F}}{\pt \hat{w}_i} = R(\th)^TW_i.  \lb{P-04}\eea
Writing $\dot{R}(\th(t))=\frac{d}{dt}R(\th(t))$, and noting that $R(\th)^T=R(\th)^{-1}$ and
$\dot{R}(\th)=\dot{\th}JR(\th)$ we obtain the function
\bea  \hat{F}_t
&\equiv& \frac{\pt \hat{F}}{\pt t} = -Z\cdot\dot{R}(\th)\hat{z} - \sum_{i=1}^{n-2}W_i\cdot\dot{R}(\th)\hat{w}_i   \nn\\
&=& -\hat{Z}\cdot R(\th)^T\dot{R}(\th)\hat{z} - \sum_{i=1}^{n-2}\hat{W}_i\cdot R(\th)^T\dot{R}(\th)\hat{w}_i   \nn\\
&=& -\dot{\th}\left(\hat{Z}\cdot R(\th)^TJR(\th)\hat{z} + \sum_{i=1}^{n-2}\hat{W}_i\cdot R(\th)^TJ{R}(\th)\hat{w}_i\right)   \nn\\
&=& -\dot{\th}\left(\hat{Z}\cdot J\hat{z} + \sum_{i=1}^{n-2}\hat{W}_i\cdot J\hat{w}_i\right).     \nn\eea
Because by the definitions (\ref{Aa}) of $A_i$s and (\ref{Bb}) of $B_{ki}$s, we obtain
\be  A_iR(\th) = R(\th)A_i, \quad B_{ki}R(\th) = R(\th)B_{ki}, \qquad \forall\,1\le i\le n-2.   \lb{P-05}\ee
By (\ref{U}), this then implies
\bea  U(z,w_1,\ldots,w_{n-2})
&=& \sum_{1\le i<j\le 3}\frac{m_im_j}{|(A_i-A_j)z+\sum_{i=1}^{n-2}(B_i-B_j)w_i|}   \nn\\
&=& \sum_{1\le i<j\le 3}\frac{m_im_j}{|(A_i-A_j)R(\th)\hat{z}+\sum_{i=1}^{n-2}(B_i-B_j)R(\th)\hat{w}_i|}   \nn\\
&=& \sum_{1\le i<j\le 3}\frac{m_im_j}{|(A_i-A_j)\hat{z}+\sum_{i=1}^{n-2}(B_i-B_j)\hat{w}_i|}   \nn\\
&=& U(\hat{z},\hat{w}_1, \ldots, \hat{w}_{n-2}),  \lb{P-06}\eea
by the orthogonality of $R(\th)$. Because $\th=\th(t)$ depends on $t$, by adding the function
$\frac{\pt \hat{F}}{\pt t}$ to the Hamiltonian function $H$ in (\ref{HaFun}), as in Line 5 in p.272
of \cite{MS}, we obtain the Hamiltonian function $\hat{H}$ in the new coordinates:
\bea
&& \hat{H}(t,\hat{Z},\hat{W}_1,\ldots,\hat{W}_{n-2},\hat{z},\hat{w}_1,\ldots,\hat{w}_{n-2})
     = H_0(Z,W_1,\ldots,W_{n-2},z,w_1,\ldots,w_{n-2}) + \hat{F}_t   \nn\\
&&\quad = \frac{1}{2}(|\hat{Z}|^2 + \sum_{i=1}^{n-2}|\hat{W}_i|^2) + (\hat{z}\cdot J\hat{Z} + \sum_{i=1}^{n-2}\hat{w}_i\cdot J\hat{W}_i)\dot{\th}
             - U(\hat{z}, \hat{w}_1,\ldots,\hat{w}_{n-2}),  \qquad \lb{P-07}\eea
where the variables of $H_0$ are functions of $\th$, $\hat{Z}$, $\hat{W}_1$, $\ldots$, $\hat{W}_{n-2}$,
$\hat{z}$, $\hat{w}_1$, $\ldots$, $\hat{w}_{n-2}$ given by (\ref{P-03})-(\ref{P-04}).

{\bf Step 2.} {\it Dilating coordinates via the polar radius $r=|z(t)|$.}

We change the coordinates $\hat{\xi}$ to
$\td{\xi}=(\td{Z},\td{W}_1,\ldots, \td{W}_{n-2}, \td{z},\td{w}_1,\ldots, \td{w}_{n-2})$
which dilate with $r=|z(t)|$ given by (\ref{rTh}). The position coordinates are transformed by
\be  \hat{z} = r\td{z}, \quad \hat{w}_i = r\td{w}_i, \qquad \forall\;1\le i\le n-2.  \lb{P-08}\ee
It is natural to scale the momenta by $1/r$ to get $\hat{Z}=\td{Z}/r$ and $\hat{W}_i=\td{W}_i/r$. But it
turns out that the new transformation with $1\le i\le n-2$
\be  \hat{Z}=\frac{1}{r}\td{Z}+\dot{r}\td{z}, \quad \hat{W}_i=\frac{1}{r}\td{W}_i+\dot{r}\td{w}_i  \lb{P-09}\ee
makes the resulting Hamiltonian function simpler. This transformation is generated by the function
\be  \td{F}(t, \td{Z}, \td{W}_1, \ldots, \td{W}_{n-2}, \hat{z}, \hat{w}_1, \ldots, \hat{w}_{n-2})
  = \frac{1}{r}(\td{Z}\cdot\hat{z} + \sum_{i=1}^{n-2}\td{W}_i\cdot\hat{w}_i) + \frac{\dot{r}}{2r}(|\hat{z}|^2 + \sum_{i=1}^{n-2}|\hat{w}_i|^2),
                        \lb{P-10}\ee
and is given by
\bea
\td{z} = \frac{\pt \td{F}}{\pt \td{Z}} = \frac{1}{r}\hat{z},    &\;&
\hat{Z} = \frac{\pt \td{F}}{\pt \hat{z}} = \frac{1}{r}\td{Z}+\frac{\dot{r}}{r}\hat{z} = \frac{1}{r}\td{Z}+\dot{r}\td{z}, \nn\\
\td{w}_i = \frac{\pt \td{F}}{\pt \td{W}_i} = \frac{1}{r}\hat{z},   &\;&
\hat{W}_i = \frac{\pt \td{F}}{\pt \hat{w}_i} = \frac{1}{r}\td{W}_i+\frac{\dot{r}}{r}\hat{w}_i = \frac{1}{r}\td{W}_i+\dot{r}\td{w}_i, \nn\eea
with
\bea  \frac{\pt \td{F}}{\pt t}
&=& -\frac{\dot{r}}{r^2}(\td{Z}\cdot\hat{z} + \sum_{i=1}^{n-2}\td{W}_i\cdot \hat{w}_i)
       + \frac{\ddot{r}r-\dot{r}^2}{2r^2}(|\hat{z}|^2 + \sum_{i=1}^{n-2}|\hat{w}_i|^2)     \nn\\
&=& -\frac{\dot{r}}{r}(\td{Z}\cdot\td{z} + \sum_{i=1}^{n-2}\td{W}_i\cdot\td{w}_i)
       + \frac{\ddot{r}r-\dot{r}^2}{2}(|\td{z}|^2 + \sum_{i=1}^{n-2}|\td{w}_i|^2),     \lb{P-11}\eea
by (\ref{P-09}).

In this case, as in the last two lines on p.272 of \cite{MS}, the Hamiltonian function $\hat{H}$ in (\ref{P-07})
becomes the new Hamiltonian function $\td{H}$ in the new coordinates:
\bea
&& \td{H}(t,\td{Z},\td{W}_1, \ldots, \td{W}_{n-2}, \td{z}, \td{w}_1, \ldots, \td{w}_{n-2})
\equiv \hat{H}(t, \hat{Z}, \hat{W}_1,\ldots,\hat{W}_{n-2}, \hat{z}, \hat{w}_1,\ldots,\hat{w}_{n-2})) + \td{F}_t \nn\\
&&\;= \frac{1}{2r^2}(|\td{Z}|^2 + \sum_{i=1}^{n-2}|\td{W}_i|^2) + \frac{\dot{r}}{r}(\td{Z}\cdot\td{z} + \sum_{i=1}^{n-2}\td{W}_i\cdot\td{w}_i)
          + \frac{\dot{r}^2}{2}(|\td{z}|^2 + \sum_{i=1}^{n-2}|\td{w}_i|^2)    \nn\\
&&\qquad\qquad  + (\td{z}\cdot J\td{Z} + \sum_{i=1}^{n-2}\td{w}_i\cdot J\td{W}_i)\dot{\th}
                - U(r\td{z},r\td{w}_i,\ldots,r\td{w}_{n-2}) + F_t   \nn\\
&&\;= \frac{1}{2r^2}(|\td{Z}|^2 + \sum_{i=1}^{n-2}|\td{W}_i|^2) + \frac{r\ddot{r}}{2}(|\td{z}|^2 + \sum_{i=1}^{n-2}|\td{w}_i|^2) \nn\\
&&\;\qquad\qquad + \; (\td{z}\cdot J\td{Z} +\sum_{i=1}^{n-2}\td{w}_i\cdot J\td{W}_i)\dot{\th}
                - \frac{1}{r}U(\td{z},\td{w}_1, \ldots, \td{w}_{n-2}).   \lb{P-12}\eea

{\bf Step 3.} {\it Coordinates via the true anomaly $\th$ as the independent variable.}

Here we want to use the true anomaly $\th\in [0,2\pi]$ as an independent variable instead of $t\in [0,T]$ to simplify
the study. This is achieved by dividing the Hamiltonian function $\td{H}$ in (\ref{P-12}) by $\dot{\th}$. Assuming
$\dot{\th}(t)>0$ for all $t\in [0,T]$, for $\td{\xi}\in W^{1,2}(\R/(T\Z), \R^8)$ we consider the action functional
corresponding to the Hamiltonian system:
\bea f(\td{\xi})
&=& \int_0^T(\frac{1}{2}\dot{\td{\xi}}(t)\cdot J\td{\xi}(t) - \td{H}(t,\td{\xi}(t))) dt  \nn\\
&=& \int_0^{2\pi}\left(\frac{1}{2}\frac{\dot{\td{\xi}}(t(\th))}{\dot{\th}(t)}\cdot J\td{\xi}(t)
          - \frac{\td{H}(t,\td{\xi}(t(\th)))}{\dot{\th}(t)}\right) d\th  \nn\\
&=& \int_0^{2\pi}\left(\frac{1}{2}\td{\xi}'(\th)\cdot J\td{\xi}(\th) - \td{H}(\th,\td{\xi}(\th))\right) d\th.  \nn\eea
Here we used $\td{\xi}'(\th)$ to denote the derivative of $\td{\xi}(\th)$ with respect to the variable $\th$. But
in the following we shall still write $\dot{\td{\xi}}(\th)$ for the derivative with respect to $\th$ instead of
$\td{\xi}'(\th)$ for notational simplicity.

It is well known that the elliptic Kepler orbit (\ref{rTh}) satisfies
$$  r(t)^2\dot{\th}(t) = \sqrt{\mu p} = \sqrt{\mu a(1-e^2)} = \sg^2 \quad {\rm with}\;\;\sg=(\mu p)^{1/4}.  $$
Note that $a=\mu^{1/3}(T/2\pi)^{2/3}$ with $T$ being the minimal period of the orbit (\ref{rTh}), we have
$$  \sg = (\mu a(1-e^2))^{1/4} = \mu^{1/3}(\frac{T}{2\pi})^{1/6}(1-e^2)^{1/4}\;\in\; (0,\mu^{1/3}(\frac{T}{2\pi})^{1/6}] $$
depending on $e$, when the mass $\mu$ and the period $T$ are fixed. Note that similarly we have $p=\sg^4/\mu$ depends
on $e$ too. Note that the function $r$ satisfies
$$  \ddot{r} = \frac{\mu p}{r^3} - \frac{\mu}{r^2} = \mu\left(\frac{p}{r^3} - \frac{1}{r^2}\right).  $$

Therefore we get the Hamiltonian function $\td{H}$ in the new coordinates:
\bea
&& \td{H}(\th,\td{Z},\td{W}_1, \ldots, \td{W}_{n-2}, \td{z}, \td{w}_1, \ldots, \td{w}_{n-2})
     \equiv \frac{1}{\dot{\th}}\td{H}(t,\td{Z},\td{W}_1, \ldots, \td{W}_{n-2}, \td{z}, \td{w}_1, \ldots, \td{w}_{n-2})   \nn\\
&&\;= \frac{1}{2r^2(t)\dot{\th}(t)}(|\td{Z}|^2 + \sum_{i=1}^{n-2}|\td{W}_i|^2)
       + \frac{r(t)\ddot{r}(t)}{2\dot{\th}(t)}(|\td{z}|^2 + \sum_{i=1}^{n-2}|\td{w}_i|^2)  \nn\\
&&\;\quad + (\td{z}\cdot J\td{Z} + \sum_{i=1}^{n-2}\td{w}_i\cdot J\td{W}_i)
          - \frac{1}{r(t)\dot{\th}(t)}U(\td{z},\td{w}_1, \ldots, \td{w}_{n-2})  \nn\\
&&\;= \frac{1}{2\sg^2}(|\td{Z}|^2 + \sum_{i=1}^{n-2}|\td{W}_i|^2) + (\td{z}\cdot J\td{Z}
        + \sum_{i=1}^{n-2}\td{w}_i\cdot J\td{W}_i) \nn\\
&&\;\quad + \frac{\mu(p-r(\th))}{2\sg^2}(|\td{z}|^2 + \sum_{i=1}^{n-2}|\td{w}_i|^2)
           - \frac{r(\th)}{\sg^2}U(\td{z},\td{w}_1, \ldots, \td{w}_{n-2}),  \lb{P-13}\eea
where $r(\th)=p/(1+e\cos\th)$. Note that now the minimal period $T$ of the elliptic solution $\td{z}=\td{z}(\th)$
becomes $2\pi$ in the new coordinates in terms of true anomaly $\th$ as an independent variable.

{\bf Step 4.} {\it Coordinates via the dilation of $\sg=(p\mu)^{1/4}$.}

The last transformation is the dilation
\be  (\td{Z},\td{W}_1, \ldots, \td{W}_{n-2}, \td{z}, \td{w}_1, \ldots, \td{w}_{n-2}) \;\mapsto\;
      (\sg\bar{Z},\sg\bar{W}_1, \ldots, \sg\bar{W}_{n-2}, \sg^{-1}\bar{z}, \sg^{-1}\bar{w}_1, \ldots, \sg^{-1}\bar{w}_{n-2}).
              \lb{P-14}\ee
This transformation is symplectic and independent of the true anomaly $\th$. Thus the
Hamiltonian function $\td{H}$ in (\ref{P-13}) becomes a new Hamiltonian function:
\bea
&& H(\th,\bar{Z},\bar{W}_1, \ldots,\bar{W}_{n-2}, \bar{z},\bar{w}_1, \ldots,\bar{w}_{n-2}) \equiv
    \td{H}(\th,\sg\bar{Z},\sg\bar{W}_1, \ldots, \sg\bar{W}_{n-2}, \sg^{-1}\bar{z}, \sg^{-1}\bar{w}_1, \ldots, \sg^{-1}\bar{w}_{n-2})  \nn\\
&&\;= \frac{1}{2}(|\bar{Z}|^2 + \sum_{i=1}^{n-2}|\bar{W}_i|^2) + (\bar{z}\cdot J\bar{Z} + \sum_{i=1}^{n-2}\bar{w}_i\cdot J\bar{W}_i)
      + \frac{p-r}{2p}(|\bar{z}|^2 + \sum_{i=1}^{n-2}|\bar{w}_i|^2)
      - \frac{r}{\sg}U(\bar{z},\bar{w}_1, \ldots,\bar{w}_{n-2}),  \qquad\quad  \lb{P-15}\eea
where one $\sg$ is factored out from $U(\sg^{-1}\bar{z}, \sg^{-1}\bar{w}_1, \ldots, \sg^{-1}\bar{w}_{n-2})$.

The proof is complete. \hb

\medskip

Motivated by ideas in Sections 2 and 3 of \cite{MS}, we now derive the linearized Hamiltonian system at such
an EEM solution of $n$-bodies, where $\sg=(\mu p)^{-1/4}$ is important.

\begin{theorem}\label{linearized.Hamiltonian}
Using notations in (\ref{PQYX}), the EEM solution $(P(t),Q(t))^T$ in time $t$ of the system (\ref{1.2}) with
\begin{equation}
Q(t)=(r(t)R(\theta(t))a_1,r(t)R(\theta(t))a_2,\ldots,r(t)R(\theta(t))a_n)^T,\quad P(t)=M\dot{Q}(t),
\end{equation}
where we denote by $M = \diag(m_1,m_1,\ldots,m_n,m_n)$,
is transformed to the new solution $(Y(\theta),X(\theta))^T$ in the true anomaly $\theta$ as the new
variable with $G=g=0$ for the original Hamiltonian function $H$ of (\ref{P21-2}), which is
given by
\be
Y(\theta)=\left(\matrix{\bar{Z}(\theta)\cr
                        \bar{W}_1(\theta)\cr
                        \ldots\cr
                        \bar{W}_{n-2}(\theta)}\right)
=\left(\matrix{0\cr
               \sigma\cr
               \ldots\cr
               \ldots\cr
               0\cr
               0}\right), \qquad
X(\theta)=\left(\matrix{\bar{z}(\theta)\cr
                        \bar{w}_1(\theta)\cr
                        \ldots\cr
                        \bar{w}_{n-2}(\theta)}\right)
=\left(\matrix{\sigma\cr
               0\cr
               \ldots\cr
               \ldots\cr
               0\cr
               0}\right).   \lb{sigma-solution}\ee

Moreover, the linearized Hamiltonian system at the EEM solution
$$   \xi_0 \equiv (Y(\theta),X(\theta))^T =
(\underbrace{0,\sigma,\ldots,\ldots,0,0}_{2(n-1)},\underbrace{\sigma,0,\ldots,\ldots,0,0}_{2(n-1)})^T\in\R^{4(n-1)}  $$
depending on the true anomaly $\theta$ with respect to the Hamiltonian function $H$ of (\ref{P21-2}) is given by
\be    \dot\zeta(\theta) = JB(\theta)\zeta(\theta),  \lb{LinearHam1}\ee
with
\bea B(\theta)
&=& H''(\theta,\bar{Z},\bar{W_1},\ldots,\bar{W}_{n-2},\bar{z},\bar{w_1},\ldots,\bar{w}_{n-2})|_{\bar\xi=\xi_0} \nn\\
&=& \left(\begin{array}{cccc|cccc}
I_2   &O     &\ldots&O     &-J    &O     &\ldots&O \\
O     &I_2   &\ldots&O     &O     &-J    &\ldots&O \\
\ldots&\ldots&\ldots&\ldots&\ldots&\ldots&\ldots&\ldots \\
O     &\ldots&O& I_2&O     &\ldots&O     &-J \\
\hline
J     &O     &\ldots&O     &H_{\bar{z}\bar{z}}(\theta,\xi_0)&O         &\ldots&O\\
O     &J     &\ldots&O     &O     &H_{\bar{w_1}\bar{w_1}}(\theta,\xi_0)&\ldots&O\\
\ldots&\ldots&\ldots&\ldots&\ldots&\ldots&\ldots&\ldots\\
O     &\ldots&O     &J     &O     &\ldots&O     &H_{\bar{w}_{n-2}\bar{w}_{n-2}}(\theta,\xi_0)
\end{array}\right),  \lb{LinearHam2}\eea
and
\be
H_{\bar{z}\bar{z}}(\theta,\xi_0) = \left(\matrix{
           -\frac{2-e\cos\theta}{1+e\cos\theta} & 0\cr
           0 & 1\cr}\right),
\quad
H_{\bar{w_i}\bar{w_i}}(\theta,\xi_0) = \left(\matrix{
           -\frac{2\bb_i+2-e\cos\theta}{1+e\cos\theta} & 0\cr
           0 & \frac{\bb_i+1+e\cos\theta}{1+e\cos\theta}\cr}\right),\ \ {\it for}\;1\le i\le n-2,  \lb{EEM}
\ee
where each $\bb_i$ with $1\le i\le n-2$ is given by (\ref{1.13}),
and $H''$ is the Hessian Matrix of $H$ with respect to its variables $\bar{Z}$,
$\bar{W_1},\ldots,\bar{W}_{n-2}$, $\bar{z}$, $\bar{w_1}$, $\ldots$, $\bar{w}_{n-2}$.
The corresponding quadratic Hamiltonian function is given by
\bea
&& H_2(\theta,\bar{Z},\bar{W_1},\ldots,\bar{W}_{n-2},\bar{z},\bar{w_1},\ldots,\bar{w}_{n-2})  \nn\\
&&\qquad \frac{1}{2}|\bar{Z}|^2+\bar{Z}\cdot J\bar{z}+\frac{1}{2}H_{\bar{z}\bar{z}}(\theta,\xi_0)|\bar{z}|^2
  + \sum_{i=1}^{n-2}\left(\frac{1}{2}|\bar{W_i}|^2+\bar{W_i}\cdot J\bar{w_i}
              +\frac{1}{2}H_{\bar{w_i}\bar{w_i}}(\theta,\xi_0)|\bar{w_i}|^2\right). \eea
\end{theorem}

{\bf Proof.} In this proof, we generalize the computations in \cite{ZL} for the EEM of the $3$-body case
to the $n$-body case here. For reader's conveniences, we given all details here. We only need to compute
$H_{\bar{z}\bar{z}}(\theta,\xi_0)$, $H_{\bar{z}\bar{w_i}}(\theta,\xi_0)$ and
$H_{\bar{w_i}\bar{w_j}}(\theta,\xi_0)$ for $1\le i,j\le n-2$ respectively.

In this proof we omit all the upper bars on the variables of $H$ in (\ref{P21-2}). By (\ref{P21-2}), we have
\bea
H_z&=&JZ+\frac{p-r}{p}z-\frac{r}{\sigma}U_z(z,w_1,\ldots,w_{n-2}),  \nn\\
H_{w_i}&=&JW_i+\frac{p-r}{p}w_i-\frac{r}{\sigma}U_{w_i}(z,w_1,\ldots,w_{n-2}),  \nn
\eea
and
\be\lb{Hessian}\left\{
\begin{array}{l}
H_{zz}=\frac{p-r}{p}I-\frac{r}{\sigma}U_{zz}(z,w_1,\ldots,w_{n-2}), \\
H_{zw_l}=H_{w_lz}=-\frac{r}{\sigma}U_{zw_l}(z,w_1,\ldots,w_{n-2}),\quad {\rm for}\;l=1,\ldots,n-2, \\
H_{w_lw_l}=\frac{p-r}{p}I-\frac{r}{\sigma}U_{w_lw_l}(z,w_1,\ldots,w_{n-2}),\quad {\rm for}\;i=l,\ldots,n-2, \\
H_{w_lw_s}=H_{w_sw_l}=-\frac{r}{\sigma}U_{w_lw_s}(z,w_1,\ldots,w_{n-2}),\quad {\rm for}\;l,s=1,\ldots,n-2,\;l\ne s,
\end{array}\right. \ee
where all the items above are $2\times2$ matrices, and we denote by $H_x$ and $H_{xy}$
the derivative of $H$ with respect to $x$, and the second derivative of $H$ with respect to
$x$ and then $y$ respectively for $x$ and $y\in\R$.

By (\ref{U_ij}) for $U_{ij}$ with $1\le i<j\le n$ and $1\le l\le n-2$, we obtain
\bea
\frac{\partial U_{ij}}{\partial z}(z,w_1,\ldots,w_{n-2})
&=& -\frac{m_im_j(a_{ix}-a_{jx})}{|(a_{ix}-a_{jx})z+\sum_{k=3}^{n}(b_{ik}-b_{jk})w_{k-2}|^3}
\left[(a_{ix}-a_{jx})z+\sum_{k=3}^{n}(b_{ik}-b_{jk})w_{k-2}\right],  \nn\\
\frac{\partial U_{ij}}{\partial w_l}(z,w_1,\ldots,w_{n-2})
&=& -\frac{m_im_j(b_{i,l+2}-b_{j,l+2})}{|(a_{ix}-a_{jx})z+\sum_{k=3}^{n}(b_{ik}-b_{jk})w_{k-2}|^3}
\left[(a_{ix}-a_{jx})z+\sum_{k=3}^{n}(b_{ik}-b_{jk})w_{k-2}\right],  \nn\eea
and
\bea
&&\frac{\partial^2 U_{ij}}{\partial z^2}(z,w_1,\ldots,w_{n-2})
=-\frac{m_im_j(a_{ix}-a_{jx})^2}{|(a_{ix}-a_{jx})z+\sum_{k=3}^{n}(b_{ik}-b_{jk})w_{k-2}|^3}I  \nn\\
&&\quad+3\frac{m_im_j(a_{ix}-a_{jx})^2}{|(a_{ix}-a_{jx})z+\sum_{k=3}^{n}(b_{ik}-b_{jk})w_{k-2}|^5}
\nn\\
&&\quad\cdot\left[(a_{ix}-a_{jx})z+\sum_{k=3}^{n}(b_{ik}-b_{jk})w_{k-2}\right]
\left[(a_{ix}-a_{jx})z+\sum_{k=3}^{n}(b_{ik}-b_{jk})w_{k-2}\right]^T,  \nn\\
&&\frac{\partial^2 U_{ij}}{\partial z\partial w_l}(z,w_1,\ldots,w_{n-2})=
-\frac{m_im_j(a_{ix}-a_{jx})(b_{i,l+2}-b_{j,l+2})}{|(a_{ix}-a_{jx})z+\sum_{k=3}^{n}(b_{ik}-b_{jk})w_{k-2}|^3}I \nn\\
&&\quad+3\frac{m_im_j(a_{ix}-a_{jx})(b_{i,l+2}-b_{j,l+2})}{|(a_{ix}-a_{jx})z+\sum_{k=3}^{n}(b_{ik}-b_{jk})w_{k-2}|^5}
\nn\\
&&\quad\cdot\left[(a_{ix}-a_{jx})z+\sum_{k=3}^{n}(b_{ik}-b_{jk})w_{k-2}\right]
\left[(a_{ix}-a_{jx})z+\sum_{k=3}^{n}(b_{ik}-b_{jk})w_{k-2}\right]^T,  \nn\\
&&\frac{\partial^2 U_{ij}}{\partial {w_l}^2}(z,w_1,\ldots,w_{n-2})
=-\frac{m_im_j(b_{i,l+2}-b_{j,l+2})^2}{|(a_{ix}-a_{jx})z+\sum_{k=3}^{n}(b_{ik}-b_{jk})w_{k-2}|^3}I  \nn\\
&&\quad+3\frac{m_im_j(b_{i,l+2}-b_{j,l+2})^2}{|(a_{ix}-a_{jx})z+\sum_{k=3}^{n}(b_{ik}-b_{jk})w_{k-2}|^5}
\nn\\
&&\quad\cdot\left[(a_{ix}-a_{jx})z+\sum_{k=3}^{n}(b_{ik}-b_{jk})w_{k-2}\right]
\left[(a_{ix}-a_{jx})z+\sum_{k=3}^{n}(b_{ik}-b_{jk})w_{k-2}\right]^T.  \nn\eea

Set
$$ K=\left(\matrix{2 & 0\cr
                   0 & -1}\right), \qquad
K_1=\left(\matrix{1 & 0\cr
                  0 & 0}\right).  $$
Now evaluating the corresponding functions at the special solution
$(\underbrace{0,\sg,\ldots,0,0}_{2(n-1)},\underbrace{\sg,0,\ldots,0,0}_{2(n-1)})^T\in\R^{4(n-1)}$
of (\ref{sigma-solution}) with $z=(\sg,0)^T$, $w_l=(0,0)^T$ for $1\le l\le n-2$, and summing them up, we
obtain
\begin{eqnarray}
\frac{\partial^2 U}{\partial z^2}\left|_{\xi_0}\right.&=&
\sum_{1\le i<j\le n}\frac{\partial^2 U_{ij}}{\partial z^2}\left|_{\xi_0}\right.
\nonumber\\
&=&\sum_{1\le i<j\le n}\left(-\frac{m_im_j(a_{ix}-a_{jx})^2}{|(a_{ix}-a_{jx})\sigma|^3}I
                        +3\frac{m_im_j(a_{ix}-a_{jx})^2}{|(a_{ix}-a_{jx})\sigma|^5}(a_{ix}-a_{jx})^2\sigma^2K_1\right)
\nonumber\\
&=&\frac{1}{\sigma^3}\left(\sum_{1\le i<j\le4}\frac{m_im_j}{|a_{ix}-a_{jx}|}\right)K
\nonumber\\
&=&\frac{\mu}{\sigma^3}K,  \label{U_zz}
\\
\frac{\partial^2 U}{\partial w_l^2}\left|_{\xi_0}\right.&=&
\sum_{1\le i<j\le n}\frac{\partial^2 U_{ij}}{\partial w_l^2}\left|_{\xi_0}\right.
\nonumber\\
&=&\sum_{1\le i<j\le n}\left(-\frac{m_im_j(b_{i,l+2}-b_{j,l+2})^2}{|(a_{ix}-a_{jx})\sigma|^3}I
                        +3\frac{m_im_j(b_{i,l+2}-b_{j,l+2})^2}{|(a_{ix}-a_{jx})\sigma|^5}(a_{ix}-a_{jx})^2\sigma^2K_1\right)
\nonumber\\
&=&\frac{1}{\sigma^3}\left(\sum_{1\le i<j\le n}\frac{m_im_j(b_{i,l+2}-b_{j,l+2})^2}{|a_{ix}-a_{jx}|^3}\right)K
\nonumber\\
&=&\frac{1}{\sigma^3}\left(\sum_{1\le i<j\le n}\frac{m_im_j(b_{i,l+2}-b_{j,l+2})b_{i,l+2}}{|a_{ix}-a_{jx}|^3}
-\sum_{1\le i<j\le n}\frac{m_im_j(b_{i,l+2}-b_{j,l+2})b_{j,l+2}}{|a_{ix}-a_{jx}|^3}\right)K
\nonumber\\
&=&\frac{1}{\sigma^3}\left(\sum_{1\le i<j\le n}\frac{m_im_j(b_{i,l+2}-b_{j,l+2})b_{i,l+2}}{|a_{ix}-a_{jx}|^3}
+\sum_{1\le j<i\le n}\frac{m_im_j(b_{i,l+2}-b_{j,l+2})b_{i,l+2}}{|a_{ix}-a_{jx}|^3}\right)K
\nonumber\\
&=&\frac{1}{\sigma^3}\left(\sum_{i=1}^n b_{i,l+2}\sum_{j=1,j\ne i}^n\frac{m_im_j(b_{i,l+2}-b_{j,l+2})}{|a_{ix}-a_{jx}|^3}\right)K
\nonumber
\\
&=&\frac{1}{\sigma^3}\left(\sum_{i=1}^n b_{i,l+2}F_{i,l+2}\right)K
\nonumber\\
&=&\frac{\mu(1+\bb_l)}{\sigma^3}K, \label{U_ww}
\end{eqnarray}
where the last equality of the first formula follows from (\ref{1.8}), and the last equality of the
second formula follows from the definition (\ref{sum.of.Fi.bi}). Similarly, we have
\begin{eqnarray}
\frac{\partial^2 U}{\partial z\partial w_l}\left|_{\xi_0}\right.
&=&\sum_{1\le i<j\le n}\frac{\partial^2 U_{ij}}{\partial z\partial w_l}\left|_{\xi_0}\right.
\nonumber\\
&=&\sum_{1\le i<j\le4}\left(-\frac{m_im_j(a_{ix}-a_{jx})(b_{i,l+2}-b_{j,l+2})}{|(a_{ix}-a_{jx})\sigma|^3}I  \right.\nn\\
& & \qquad\qquad \left.+3\frac{m_im_j(a_{ix}-a_{jx})(b_{i,l+2}-b_{j,l+2})}{|(a_{ix}-a_{jx})\sigma|^5}(a_{ix}-a_{jx})^2\sigma^2K_1\right)
\nonumber\\
&=&\left(\sum_{1\le i<j\le n}\frac{m_im_j(b_{i,l+2}-b_{j,l+2})\cdot sign(a_{ix}-a_{jx})}{|(a_{ix}-a_{jx})|^2}\right)\frac{K}{\sigma^3}
\nonumber\\
&=&\left(\sum_{1\le i<j\le n}\frac{m_im_j(b_{i,l+2}-b_{j,l+2})\cdot(-1)}{|(a_{ix}-a_{jx})|^2}\right)\frac{K}{\sigma^3}
\nonumber\\
&=&\left(-\sum_{1\le i<j\le n}\frac{m_im_jb_{i,l+2}}{|(a_{ix}-a_{jx})|^2}
+\sum_{1\le i<j\le n}\frac{m_im_jb_{j,l+2}}{|(a_{ix}-a_{jx})|^2}\right)\frac{K}{\sigma^3}
\nonumber\\
&=&\left(-\sum_{i=1}^nb_{i,l+2}\sum_{j=i+1}^n\frac{m_im_j}{|(a_{ix}-a_{jx})|^2}
+\sum_{j=1}^nb_{j,l+2}\sum_{i=1}^{j-1}\frac{m_im_j}{|(a_{ix}-a_{jx})|^2}\right)\frac{K}{\sigma^3}
\nonumber\\
&=&\left(-\sum_{i=1}^nb_{i,l+2}\sum_{j=i+1}^n\frac{m_im_j}{|(a_{ix}-a_{jx})|^2}
+\sum_{i=1}^nb_{i,l+2}\sum_{j=1}^{i-1}\frac{m_im_j}{|(a_{ix}-a_{jx})|^2}\right)\frac{K}{\sigma^3}
\nonumber\\
&=&\left(\sum_{i=1}^nb_{i,l+2}\sum_{j=i+1}^n\frac{m_im_j(a_{ix}-a_{jx})}{|(a_{ix}-a_{jx})|^3}
+\sum_{i=1}^nb_{i,l+2}\sum_{j=1}^{i-1}\frac{m_im_j(a_{ix}-a_{jx})}{|(a_{ix}-a_{jx})|^3}\right)\frac{K}{\sigma^3}
\nonumber\\
&=&-\left(\sum_{i=1}^nb_{i,l+2}\sum_{j=1,j\ne i}^n\frac{m_im_j(a_{ix}-a_{jx})}{|(a_{ix}-a_{jx})|^3}\right)\frac{K}{\sigma^3}
\nonumber\\
&=&-\left(\sum_{i=1}^nb_{i,l+2}\mu m_ia_{ix}\right)\frac{K}{\sigma^3}
\nonumber\\
&=&O,\label{U_zw}
\end{eqnarray}
where in the fourth and fourth last equality, we used the ascending order of $a_{ix},1\le i\le n$,
in the second last equation, we used (\ref{eq.of.cc}), and in the last equality, we used (\ref{vMv}).
Moreover, for $l\ne s$, we have
\begin{eqnarray}
\frac{\partial^2 U}{\partial w_l\partial w_s}\left|_{\xi_0}\right.&=&
\sum_{1\le i<j\le n}\frac{\partial^2 U_{ij}}{\partial w_l\partial w_s}\left|_{\xi_0}\right.
\nonumber\\
&=&\sum_{1\le i<j\le n}\left(-\frac{m_im_j(b_{i,l+2}-b_{j,l+2})(b_{i,s+2}-b_{j,s+2})}{|(a_{ix}-a_{jx})\sigma|^3}I\right.
\nonumber
\\
&&+\left.3\frac{m_im_j(b_{i,l+2}-b_{j,l+2})(b_{i,s+2}-b_{j,s+2})}{|(a_{ix}-a_{jx})\sigma|^5}(a_{ix}-a_{jx})^2\sigma^2K_1\right)
\nonumber\\
&=&\frac{1}{\sigma^3}\left(\sum_{1\le i<j\le n}\frac{m_im_j(b_{i,l+2}-b_{j,l+2})(b_{i,s+2}-b_{j,s+2})}{|a_{ix}-a_{jx}|^3}\right)K
\nonumber\\
&=&\frac{1}{\sigma^3}\left(\sum_{1\le i<j\le n}\frac{m_im_jb_{i,l+2}(b_{i,s+2}-b_{j,s+2})}{|a_{ix}-a_{jx}|^3}
-\sum_{1\le i<j\le n}\frac{m_im_jb_{j,l+2}(b_{i,s+2}-b_{j,s+2})}{|a_{ix}-a_{jx}|^3}\right)K
\nonumber\\
&=&\frac{1}{\sigma^3}\left(\sum_{1\le i<j\le n}\frac{m_im_jb_{i,l+2}(b_{i,s+2}-b_{j,s+2})}{|a_{ix}-a_{jx}|^3}
+\sum_{1\le j<i\le n}\frac{m_jm_ib_{i,l+2}(b_{i,s+2}-b_{j,s+2})}{|a_{jx}-a_{ix}|^3}\right)K
\nonumber\\
&=&\frac{1}{\sigma^3}\left(\sum_{i=1}^nb_{i,l+2}\sum_{j=1,j\ne i}^n\frac{m_im_j(b_{i,s+2}-b_{j,s+2})}{|a_{ix}-a_{jx}|^3}\right)K
\nonumber
\\
&=&\frac{1}{\sigma^3}(\sum_{i=1}^nb_{i,l+2}F_{i,s+2})K
\nonumber\\
&=&\frac{1}{\sigma^3}(\sum_{i=1}^nb_{i,l+2}(\mu-\lambda_{s+2})m_ib_{i,s+2})K
\nonumber\\
&=&O,\label{U_w1w2}
\end{eqnarray}
where in the third last equality, we used (\ref{F_ki}),
and in the last equality of (\ref{U_w1w2}), we used (\ref{vMv}) and (\ref{bi}).

By (\ref{U_zz}), (\ref{U_ww}), (\ref{U_zw}) and (\ref{Hessian}), we have
\begin{eqnarray}
H_{zz}|_{\xi_0}&=&\frac{p-r}{p}I-\frac{r\mu}{\sigma^4}K=I-\frac{r}{p}I-\frac{r\mu}{p\mu}K
=I-\frac{r}{p}(I+K)
=\left(\matrix{-\frac{2-e\cos\theta}{1+e\cos\theta} & 0\cr
               0 & 1}\right),  \nn\\
H_{zw_l}|_{\xi_0}&=&-\frac{r}{\sigma}\frac{\partial^2U}{\partial z\partial w_l}|_{\xi_0}=O,\quad {\rm for}\;1\le l\le n-2,
\nn\\
H_{w_lw_s}|_{\xi_0}&=&-\frac{r}{\sigma}\frac{\partial^2U}{\partial w_l\partial w_s}|_{\xi_0}=O,
\quad {\rm for}\;1\le l,s\le n-2,\;l\ne s,\nn\\
H_{w_lw_l}|_{\xi_0}&=&\frac{p-r}{p}I-\frac{r(1+\bb_l)\mu}{\sigma^4}K
=I-\frac{r}{p}I-\frac{r(1+\bb_l)\mu}{p\mu}K
\nn\\
&=&I-\frac{r}{p}(I+(1+\bb_l) K)=
\left(\matrix{-\frac{2\bb_l+2-e\cos\theta}{1+e\cos\theta} & 0\cr
              0 & \frac{\bb_l+1+e\cos\theta}{1+e\cos\theta}}\right),
\quad {\rm for}\;1\le l\le n-2.
\end{eqnarray}
Thus the proof is complete.\hb

\begin{remark}\label{remark:bb1}
(i) When we set $n=3$ in Theorem \ref{linearized.Hamiltonian}, then $\bb_1$ is precisely the mass
parameter $\beta$ defined by (\ref{1.4}), and the corresponding linearized Hamiltonian system at
the EEM $q_{m,e}(t)$ is given by
\be
z'=J\left(\matrix{1 & 0 & 0 & 1\cr
                       0 & 1 & -1 & 0 \cr
                       0 & -1 &\frac{-2\beta-2+e\cos(t)}{1+e\cos(t)} & 0 \cr
                       1 & 0 & 0 & \frac{\beta+1+e\cos(t)}{1+e\cos(t)} \cr}\right)z. \label{2.41}
\ee
Note that this system was derived in \cite{MSS1} and \cite{ZL} too.

(ii) The Hamiltonian equation of the $i$-th part of the other $(n-2)$ parts with $1\le i\le n-2$
is given by
\be
z'=J\left(\matrix{1 & 0 & 0 & 1\cr
                       0 & 1 & -1 & 0 \cr
                       0 & -1 &\frac{-2\beta_i-2+e\cos(t)}{1+e\cos(t)} & 0 \cr
                       1 & 0 & 0 & \frac{\beta_i+1+e\cos(t)}{1+e\cos(t)} \cr}\right)z. \label{2.42}
\ee
Also, $\bb_1$ coincides with $\bb_c$ in Table 2 of \cite{MSS1} when $\alpha=1$.
\end{remark}

Now we can give

{\bf Proof of Theorem \ref{T1.2}.} Note that by Theorem \ref{linearized.Hamiltonian}, specially
(\ref{LinearHam1})-(\ref{EEM}), we obtain that the matrix $H_{\bar{z}\bar{z}}(\th,\xi_0)$ together
with the first identity matrix $I_2$ in the diagonal of the matrix $B(\th)$ in (\ref{LinearHam2})
yield a $4$-dimensional Hamiltonian system corresponding to the Kepler $2$-body problem, and each
matrix $H_{\bar{w}_i\bar{w}_i}(\th,\xi_0)$ together with the $(i+1)$-th identity matrix $I_2$ in
the diagonal of the matrix $B(\th)$ in (\ref{LinearHam2}) yield a $4$-dimensional Hamiltonian system
(\ref{2.42}) with $\beta_i$ given by (\ref{1.13}), which corresponds to the linear system (\ref{2.41})
of the Euler $3$-body problem with $\beta$ replaced by $\beta_i$ for $1\le i\le n-2$. Therefore
Theorem \ref{T1.2} holds. \hb

\setcounter{equation}{0}
\section{A collinear $4$-body problem with two small masses in the middle}\label{sec:3}

We now consider the linear stability of special collinear central configurations in the four
body problem with two small masses in the middle. A typical example is the EEM orbit of the
$4$-bodies, the Earth, the Moon and two space stations in the middle as mentioned at the
beginning of this paper with $n=4$. We try to give an analytical way following which one can
numerically find out the best elliptic-hyperbolic positions for the two space stations using
results in \cite{MSS1} and \cite{MSS2}. Specially, for the four masses we fix $m_1=m\in (0,1)$,
and let $m_2=\ep$, $m_3=\tau\ep$, $m_4=1-m-(\tau+1)\ep$ with $\tau>0$ and
$0<\ep<\frac{1-m}{\tau+1}$. They satisfy
\be  m_1+m_2+m_3+m_4=1.  \lb{3.1}\ee
Suppose $q_1$, $q_2$, $q_3$ and $q_4$ are four points on the $x$-axis in $\R^2$, and
form a central configuration. Using notations similar to those in \cite{ZL}, we set
\be  q_1=(0,0)^T,\quad q_2=(x\aa,0)^T,\quad q_3=(y\aa,0)^T,\quad q_4=(\aa,0)^T, \lb{3.2}\ee
where $\aa=\aa_{\ep,\tau}=|q_4-q_1|$, $x=x_{\ep,\tau}$, $y=y_{\ep,\tau}$ satisfy
$0<x<y<1$. Then the center of mass of the four particles is
\bea  q_c
&=& m_1q_1+m_2q_2+m_3q_3+m_4q_4   \nn\\
&=& ([m_2x+m_3y+m_4]\aa,0)^T    \nn\\
&=& ([(1-m)-(1+\tau-x-\tau y)\ep]\aa,0)^T,  \lb{3.3}\eea
where (\ref{3.1}) is used to get the last equality.

For $i=1$, $2$, $3$ and $4$, let $a_i=q_i-q_c$, and denote by $a_{ix}$ and $a_{iy}$
the $x$ and $y$-coordinates of $a_i$ respectively. Then we have
\begin{eqnarray}\label{a_ix}
a_{1x} &=& [m-1+(1+\tau-x-\tau y)\epsilon]\alpha,\quad\quad\quad\;\;
            a_{2x}=[m+x-1+(1+\tau -x-\tau y)\epsilon]\alpha, \lb{3.4}\\
a_{3x} &=& [m+y-1+(1+\tau -x-\tau y)\epsilon]\alpha,\quad\quad
            a_{4x}=[m+(1+\tau -x-\tau y)\epsilon]\alpha,   \lb{3.5}\end{eqnarray}
and
\be  a_{iy}=0,\quad\quad {\rm for}\; i=1,2,3,4.    \lb{3.6}\ee
Next we study properties of this central configuration.

\medskip

{\bf Step 1.} {\it Computations on $\aa$ and $x,y$.}

Scaling $\alpha$ by setting $\sum_{i=1}^4m_i|a_i|^2=1$, we obtain
\begin{eqnarray}  \frac{1}{\alpha^2}
&=& \frac{\sum_{i=1}^4m_i|a_i|^2}{\alpha^2}   \nonumber\\
&=& m[m-1+(1+\tau -x-\tau y)\epsilon]^2+\epsilon[m+x-1+(1+\tau -x-\tau y)\epsilon]^2  \nonumber\\
& &\quad +\tau \epsilon[m+y-1+(1+\tau -x-\tau y)\epsilon]^2+(1-m-(\tau +1)\epsilon)[m+(1+\tau -x-\tau y)\epsilon]^2  \nonumber\\
&=& m(1-m)^2+2m(1-m)(x+\tau y-1-\tau )\epsilon+m(x+\tau y-1-\tau )^2\epsilon^2  \nonumber\\
& &\quad +\epsilon[(x-1)^2+(m+(1+\tau -x-\tau y)\epsilon)^2+2(x-1)(m+(1+\tau -x-\tau y)\epsilon)]  \nonumber\\
& &\quad +\tau \epsilon[(y-1)^2+(m+(1+\tau -x-\tau y)\epsilon)^2+2(y-1)(m+(1+\tau -x-\tau y)\epsilon)]  \nonumber\\
& &\quad +(1-m)[m^2+2m(1+\tau -x-\tau y)+(1+\tau -x-\tau y)^2\epsilon^2]-(\tau +1)\epsilon(m+(1+\tau -x-\tau y)\epsilon)^2  \nonumber\\
&=& m(1-m) + (1+\tau -x-\tau y)^2\epsilon^2 + \epsilon[(x-1)^2+\tau (y-1)^2]  \nn\\
& &\quad  + \ep m(1+\tau )(m+(1+\tau -x-\tau y)\ep)^2 - 2\ep(1+\tau -x-\tau y)(m+(1+\tau -x-\tau y)\ep)  \nn\\
&=& m(1-m) + [(1-x)^2 + \tau (1-y)^2 + m^3(1+\tau ) - 2m(1+\tau -x-\tau y)]\ep  \nn\\
& &\quad + [2m^2(1+\tau )(1+\tau -x-\tau y) - (1+\tau -x-\tau y)^2]\ep^2 + m(1+\tau )(1+\tau -x-\tau y)^2\ep^3.   \lb{3.7}\eea
Moreover, let
\be  \alpha_0=\lim_{\epsilon\to 0}\alpha=[m(1-m)]^{-{1\over2}},  \lb{3.8}\ee
and
\be  q_{c,0}=\lim_{\epsilon\to 0}q_c=(1-m)\alpha_0,   \lb{3.9}\ee
and hence
\bea
a_{1x,0} &=& \lim_{\epsilon\to 0}a_{1x} = -(1-m)\alpha_0,  \lb{3.10}\\ 
a_{4x,0} &=& \lim_{\epsilon\to 0}a_{4x} = m\alpha_0.   \lb{3.11}\eea 
The potential $\mu$ is given by
\be
\mu=\mu_{\epsilon,\tau}=\sum_{1\le i<j\le4}\frac{m_im_j}{|a_i-a_j|},
\ee
and by Lemma 3 of \cite{IM}, we have
\be
\mu_0=\lim_{\epsilon\to 0}\mu=\frac{m(1-m)}{\alpha_0}=\alpha_0^{-3}.\label{mu0}
\ee
In the following, we will use the subscript $0$ to denote the limit value of the parameters
when $\epsilon\to 0$.

Motivated by Proposition 1 in \cite{Xia}, we have

\medskip

\begin{lemma}\label{Lemma3.1}
When $\epsilon\rightarrow0$, $a_2$ and $a_3$ must converge to the same point $a^*$.
Moreover, $a_{1,0},a^*,a_{4,0}$ is the central configuration of the restricted 3-body problem
with given masses $\tilde{m}_1=m,\tilde{m}_2=0,\tilde{m}_3=1-m$
which the small mass lies in the segment between the other two masses.
\end{lemma}

\medskip

{\bf Proof.} If $a_2$ and $a_3$ do not converge to the same point when $\ep\to 0$, there is a
sequence $\{\epsilon_n\}_{n=0}^\infty$ convergent to $0$ such that
\be   a_{23x}^*\equiv\lim_{n\to \infty}(a_{2x}-a_{3x}) \ne 0.  \lb{notZero}\ee
Up to a subsequence of $\{\epsilon_n\}_{n=0}^\infty$, and we denote it still by
$\{\epsilon_n\}_{n=0}^\infty$, we have
\be
\lim_{n\to \infty}a_{2x}=a_{2x}^*.
\ee
Then
\be
\lim_{n\to \infty}a_{3x}=\lim_{n\to \infty}a_{2x}+\lim_{n\to \infty}(a_{2x}-a_{3x})=a_{2x}^*+a_{23x}^*.
\ee

Because $a_1$, $a_2$, $a_3$ and $a_4$ form a central configuration, for the two middle points we have
\bea
\frac{m_1(a_2-a_1)}{|a_2-a_1|^3}+\frac{m_3(a_2-a_3)}{|a_2-a_3|^3}+\frac{m_4(a_2-a_4)}{|a_2-a_4|^3}&=&\mu a_2,
\label{cc.eq1}
\\
\frac{m_1(a_3-a_1)}{|a_3-a_1|^3}+\frac{m_2(a_3-a_2)}{|a_3-a_2|^3}+\frac{m_4(a_3-a_4)}{|a_3-a_4|^3}&=&\mu a_3.
\label{cc.eq2}
\eea
Let $\epsilon=\epsilon_n,n\in\N$, and $n\to\infty$,
together with (\ref{3.5}), (\ref{3.9}), (\ref{3.10}), (\ref{mu0}) and (\ref{notZero}),
(\ref{cc.eq1}) and (\ref{cc.eq2}) become
\bea
\frac{m}{(a_{2x}^*-a_{1x,0})^2}-\frac{1-m}{(a_{2x}^*-a_{4x,0})^2} &=&\mu_0 a_{2x}^*, \label{lim.cc.eq1}\\
\frac{m}{(a_{2x}^*+a_{23x}^*-a_{1x,0})^2}-\frac{1-m}{(a_{2x}^*+a_{23x}^*-a_{4x,0})^2} &=&\mu_0(a_{2x}^*+a_{23x}^*).
\label{lim.cc.eq2}
\eea

We define
\be
f(t)=\frac{m}{(t-a_{1x,0})^2}-\frac{1-m}{(t-a_{4x,0})^2}-\mu_0 t, \qquad {\rm for}\;\;t\in (a_{1x,0},a_{4x,0}).
\ee
Then $f$ is a strictly decreasing function satisfying
\be  \lim_{t\to a_{1x,0}}f(t)=+\infty, \quad {\rm and}\quad  \lim_{t\to a_{4x,0}}f(t)=-\infty.   \lb{bryCond}\ee
Thus there is a unique zero point of $f$ in $[a_{1x,0},a_{4x,0}]$, which we denote by $a_x^*$.
Here (\ref{bryCond}) yields $a_{1x,0} < a_x^* < a_{4x,0}$.

Now (\ref{lim.cc.eq1}) and (\ref{lim.cc.eq2}) yield two zero points $a_{2x}^*$ and
$a_{2x}^*+a_{23x}^*$ of $f$ in $[a_{1x,0},a_{4x,0}]$ respectively, which then yields a
contradiction. Therefore, we have proved $\lim_{\epsilon\to 0}(a_{2x}-a_{3x})=0$.

Now we want to prove $\lim_{\epsilon\to0}a_{2x}=a_x^*$.
If not, there is a sequence $\{\tilde\epsilon_n\}_{n=0}^\infty$ converges to $0$,
such that $\lim_{n\to \infty}a_{2x}=\tilde{a}^*\ne a^*$.
Then $\lim_{n\to \infty}a_{3x}=\tilde{a}^*$.

Now adding $\frac{m_2}{m_2+m_3}$ times (\ref{cc.eq1}) to $\frac{m_3}{m_2+m_3}$ times (\ref{cc.eq2})
yields
\be
m_1\left(\frac{1}{\tau+1}\frac{a_2-a_1}{|a_2-a_1|^3}+\frac{\tau}{\tau+1}\frac{a_3-a_1}{|a_3-a_1|^3}\right)
+m_4\left(\frac{1}{\tau+1}\frac{a_2-a_4}{|a_2-a_4|^3}+\frac{\tau}{\tau+1}\frac{a_3-a_4}{|a_3-a_4|^3}\right)
=\mu\frac{a_2+\tau a_3}{1+\tau},
\label{cc.eq12}
\ee
Let $\epsilon=\tilde\epsilon_n,n\in\N$, and $n\to\infty$,
together with (\ref{3.6}), (\ref{3.10}), (\ref{3.11}) and (\ref{mu0}),
(\ref{cc.eq12}) becomes
\be
\frac{m}{(\tilde{a}^*-a_{1x,0})^2}-\frac{1-m}{(\tilde{a}^*-a_{4x,0})^2}
=\mu_0 \tilde{a}^*,
\ee
then using also the property of unique zero point of $f(x)$, we obtain a contradiction.
Thus we must have $\lim_{\epsilon\to 0}a_{2x} = \lim_{\epsilon\to 0}a_{3x} = a_x^*$.

By direct computations, we can check that $a_{1,0}=(a_{1x,0},0)^T$, $a^*=(a_x^*,0)^T$ and
$a_{4,0}=(a_{4x,0},0)^T$ form a collinear central configuration with given masses
$\tilde{m}_1=m$, $\tilde{m}_2=0$ and $\tilde{m}_3=1-m$. The uniqueness is obtained by these
three given ordered masses as in \cite{Mou}. \hb

\medskip

By Lemma \ref{Lemma3.1}, we can suppose
\be
\lim_{\epsilon\rightarrow0}x=\lim_{\epsilon\rightarrow0}y=x_0,  \lb{x_0def}
\ee
and hence
\bea
a_{2x,0} &=& \lim_{\epsilon\rightarrow0}a_{2x}=(m+x_0-1)\alpha_0, \\
a_{3x,0} &=& \lim_{\epsilon\rightarrow0}a_{3x}=(m+x_0-1)\alpha_0. \eea
Note that $a_{1,0}$, $a_{2,0}$ and $a_{4,0}$ form a central configuration with given masses
$\tilde{m}_1=m$, $\tilde{m}_2=0$ and $\tilde{m}_3=1-m$. Then $\tilde{x}=\frac{x_0}{1-x_0}$ is the
unique positive root of Euler's quintic polynomial equation (cf. p. 276 of \cite{Win1} and p.29 of \cite{Lon}):
\begin{equation}\label{quintic.polynomial}
(1-m)\tilde{x}^5+(3-3m)\tilde{x}^4+(3-3m)\tilde{x}^3-3m\tilde{x}^2-3m\tilde{x}-m=0.
\end{equation}
Thus $x_0$ satisfies:
\be  x_0^5-(3-m)x_0^4+(3-2m)x_0^3-mx_0^2+2mx_0-m=0.  \lb{x0Equation}\ee

Next we derive the equations satisfied by $x=x(\ep)$ and $y=y(\ep)$. Because $a_1$, $a_2$,
$a_3$ and $a_4$ form a central configuration, we have
\bea
\frac{\epsilon}{x^2\alpha^2}+\frac{\tau\epsilon}{y^2\alpha^2}+\frac{1-m-(1+\tau)\epsilon}{\alpha^2}
  &=& \mu[1-m-(1+\tau-x-\tau{y})\epsilon]\alpha, \label{cc.eq1'}\\
-\frac{m}{x^2\alpha^2}+\frac{\tau\epsilon}{(y-x)^2\alpha^2}+\frac{1-m-(1+\tau)\epsilon}{(1-x)^2\alpha^2}
  &=& \mu[1-m-x-(1+\tau-x-\tau{y})\epsilon]\alpha, \label{cc.eq2'}\\
-\frac{m}{y^2\alpha^2}-\frac{\epsilon}{(y-x)^2\alpha^2}+\frac{1-m-(1+\tau)\epsilon}{(1-y)^2\alpha^2}
  &=& \mu[1-m-y-(1+\tau-x-\tau{y})\epsilon]\alpha, \label{cc.eq3'}\\
-\frac{m}{\alpha^2}-\frac{\epsilon}{(1-x)^2\alpha^2}-\frac{\tau\epsilon}{(1-y)^2\alpha^2}
  &=&-\mu[m+(1+\tau-x-\tau{y})\epsilon]\alpha. \label{cc.eq4'}\eea
From (\ref{cc.eq1'}) and (\ref{cc.eq2'}), we have
\bea  0
&=& \left[\frac{\epsilon}{x^2}+\frac{\tau\epsilon}{y^2}+1-m-(1+\tau)\epsilon\right][1-m-x-(1+\tau-x-\tau{y})\epsilon]
\nonumber\\
& & -\left[-\frac{m}{x^2}+\frac{\tau\epsilon}{(y-x)^2}+\frac{1-m-(1+\tau)\epsilon}{(1-x)^2}\right][1-m-(1+\tau-x-\tau{y})\epsilon]
\nonumber\\
&=& (1-m)(1-m-x)-\left[-\frac{m}{x^2}+\frac{1-m}{(1-x)^2}\right](1-m)  \nonumber\\
& & +\epsilon\left[-(1-m)(1+\tau-x-\tau{y})+(1-m-x)(\frac{1}{x^2}+\frac{\tau}{y^2}-1-\tau)\right.
\nonumber\\
& & \quad\left.-(1-m)(\frac{\tau}{(y-x)^2}-\frac{1+\tau}{(1-x)^2})+(-\frac{m}{x^2}+\frac{1-m}{(1-x)^2})(1+\tau-x-\tau{y})\right]
\nonumber\\
& & +\epsilon^2\left[(\frac{1}{x^2}+\frac{\tau}{y^2}-1-\tau)(1+\tau-x-\tau{y})+(\frac{\tau}{(y-x)^2}
           -\frac{1+\tau}{(1-x)^2})(1+\tau-x-\tau{y})\right]
\nonumber\\
&=& -(1-m)\frac{x^5-(3-m)x^4+(3-2m)x^3-mx^2+2mx-m}{x^2(1-x)^2}  \nonumber\\
& & +\epsilon\left[(m-1-\frac{m}{x^2}+\frac{1-m}{(1-x)^2})(1+\tau-x-\tau{y})+(1-m-x)(\frac{1}{x^2}+\frac{\tau}{y^2}-1-\tau)\right.
\nonumber\\
& & \quad\left.-(1-m)(\frac{\tau}{(y-x)^2}-\frac{1+\tau}{(1-x)^2})\right]
  +\epsilon^2\left[\frac{1}{x^2}+\frac{\tau}{y^2}+\frac{\tau}{(y-x)^2}-\frac{1+\tau}{(1-x)^2}-1-\tau\right](1+\tau-x-\tau{y}).
\nonumber\\
\label{g.epsilon}
\eea
We denote the right hand side of (\ref{g.epsilon}) by $g_\epsilon(x,y)$,
then $x^2(1-x)^2y^2(y-x)^2g_\epsilon(x,y)$ is a binary polynomial in $x,y$.
Similarly, from (\ref{cc.eq1'}) and (\ref{cc.eq3'}), we have
\be
h_\epsilon(x,y)=0,
\ee
where
\bea
h_\epsilon(x,y)&=&-(1-m)\frac{y^5-(3-m)y^4+(3-2m)y^3-my^2+2my-m}{y^2(1-y)^2}
\nonumber
\\
&&+\epsilon\left[(m-1-\frac{m}{y^2}+\frac{1-m}{(1-y)^2})(1+\tau-x-\tau{y})+(1-m-y)(\frac{1}{x^2}+\frac{\tau}{y^2}-1-\tau)\right.
\nonumber
\\
&&\quad\left.+(1-m)(\frac{1}{(y-x)^2}+\frac{1+\tau}{(1-y)^2})\right]
\nonumber
\\
&&+\epsilon^2\left[\frac{1}{x^2}+\frac{\tau}{y^2}-\frac{1}{(y-x)^2}-\frac{1+\tau}{(1-y)^2}-1-\tau\right](1+\tau-x-\tau{y}).
\label{h.epsilon}
\eea
Therefore, $x$ and $y$ can be solved out from $g_\epsilon(x,y)=0$ and $h_\epsilon(x,y)=0$.

Now by the first conclusion of Lemma 3.1, letting $\ep\to 0$ in the equations $g_\epsilon(x,y)=0$ and
$h_\epsilon(x,y)=0$, we obtain
\bea
&& \lim_{\epsilon\to 0}x^2(1-x)^2y^2(y-x)^2g_\epsilon(x,y)  \nonumber\\
&& \quad = -(1-m)y^2(y-x)^2[x^5-(3-m)x^4+(3-2m)x^3-mx^2+2mx-m], \label{lim.g.epsilon}\\
&& \lim_{\epsilon\to 0}y^2(1-y)^2x^2(x-y)^2h_\epsilon(x,y)  \nonumber\\
&& \quad = -(1-m)x^2(x-y)^2[y^5-(3-m)y^4+(3-2m)y^3-my^2+2my-m]. \label{lim.h.epsilon}\eea

Here in (\ref{lim.g.epsilon}) and (\ref{lim.h.epsilon}) we have the same polynomial
again as that in the left hand side of (\ref{x0Equation}).

\medskip

{\bf Step 2.} {\it Computations on $\bb_i$s.}

Now in our case, $D$ is given by
\begin{equation}
D=
\left(
\matrix{
\mu-\frac{1}{\alpha^3}[\frac{\epsilon}{x^3}+\frac{\tau\epsilon}{y^3}+1-m-(1+\tau)\epsilon],\quad\quad
         \frac{\epsilon}{x^3\alpha^3},\quad\quad
         \frac{\tau\epsilon}{y^3\alpha^3},\quad\quad \frac{1-m-(1+\tau)\epsilon}{\alpha^3}\cr
\frac{m}{x^3\alpha^3},\quad\quad
         \mu-\frac{1}{\alpha^3}[\frac{m}{x^3}+\frac{\tau\epsilon}{(y-x)^3}+\frac{1-m-(1+\tau)\epsilon}{(1-x)^3}],\quad\quad
         \frac{\tau\epsilon}{(y-x)^3\alpha^3},\quad\quad     \frac{1-m-(1+\tau)\epsilon}{(1-x)^3\alpha^3}\cr
\frac{m}{y^3\alpha^3},\quad\quad     \frac{\epsilon}{(y-x)^3\alpha^3},\quad\quad
         \mu-\frac{1}{\alpha^3}[\frac{m}{y^3}+\frac{\epsilon}{(y-x)^3}+\frac{1-m-(1+\tau)\epsilon}{(1-y)^3}],\quad\quad
         \frac{1-m-(1+\tau)\epsilon}{(1-y)^3\alpha^3}\cr
\frac{m}{\alpha^3},\quad\quad\quad           \frac{\epsilon}{(1-x)^3\alpha^3},\quad\quad\quad
         \frac{\tau\epsilon}{(1-y)^3\alpha^3},\quad\quad\quad\quad
         \mu-\frac{1}{\alpha^3}[m+\frac{\epsilon}{(1-x)^3}+\frac{\tau\epsilon}{(1-y)^3}]}\right).
\label{D_epsilon}
\end{equation}
Recall that the other two eigenvalues of $D$ are $\lambda_3$ and $\lambda_4$, then we have
\bea
\det(D-\lambda I_4)&=&-\lambda(\mu-\lambda)(\lambda_3-\lambda)(\lambda_4-\lambda)
\nonumber\\
&=&\lambda^4-(\mu+\lambda_3+\lambda_4)\lambda^3+(\lambda_3\lambda_4+\mu(\lambda_3+\lambda_4))\lambda^2-\lambda_3\lambda_4\mu\lambda
\eea
On the other hand
\bea
\det(D-\lambda I_4)&=&\lambda^4-(\emph{tr}D)\lambda^3+(\sum_{i,j=1,i<j}^4\det E_{ij})\lambda^2+\ldots,
\eea
where $E_{ij}$ is the principal minor when deleting all the rows and columns except for $i$ and $j$.
Then we have
\bea
\mu+\lambda_3+\lambda_4&=&\emph{tr}D
\nonumber\\
&=&\sum_{i=1}^4 D_{ii}
\nonumber\\
&=&\sum_{i=1}^4\left(\mu-\sum_{j=1,j\ne i}^4 D_{ij}\right)
\nonumber\\
&=&4\mu-\sum_{i,j=1,i\ne j}^4 D_{ij}
\nonumber\\
&=&4\mu-\frac{1}{\alpha^3}\left[\frac{m+\epsilon}{x^3}+\frac{m+\tau\epsilon}{y^3}+1-(1+\tau)\epsilon
  +\frac{(1+\tau)\epsilon}{(y-x)^3}+\frac{1-m-\tau\epsilon}{(1-x)^3}+\frac{1-m-\epsilon}{(1-y)^3}\right],
\nonumber
\\
\\
\lambda_3\lambda_4+\mu(\lambda_3+\lambda_4)&=&\sum_{i,j=1,i<j}^4\det E_{ij}
\nonumber\\
&=&\sum_{i,j=1,i<j}^4\left(D_{ii}D_{jj}-D_{ij}D_{ji}\right)
\nonumber\\
&=&\sum_{i,j=1,i<j}^4D_{ii}D_{jj}-\sum_{i,j=1,i<j}^4D_{ij}D_{ji}
\nonumber\\
&=&\frac{\left(\sum_{i=1}^4D_{ii}\right)^2-\sum_{i=1}^4D_{ii}^2}{2}-\sum_{i,j=1,i<j}^4D_{ij}D_{ji}
\nonumber\\
&=&{1\over2}(trD)^2-{1\over2}\sum_{i=1}^4\left(\mu-\sum_{s=1,s\ne i}^4D_{is}\right)^2
-\sum_{i,j=1,i<j}^4D_{ij}D_{ji}
\nonumber\\
&=&{1\over2}(trD)^2-2\mu^2+\mu\sum_{i,j=1,i\ne j}^4 D_{ij}
  -{1\over2}\sum_{i=1}^4\left(\sum_{s=1,s\ne i}^4D_{is}\right)^2
  -\sum_{i,j=1,i<j}^4D_{ij}D_{ji}
\eea

Let
\bea
\delta&=&\frac{1}{2\mu}\sum_{i,j=1,i\ne j}^4 D_{ij}
\nonumber
\\
&=&\frac{1}{2\mu\alpha^3}\left[\frac{m+\epsilon}{x^3}+\frac{m+\tau\epsilon}{y^3}+1-(1+\tau)\epsilon
+\frac{(1+\tau)\epsilon}{(y-x)^3}+\frac{1-m-\tau\epsilon}{(1-x)^3}+\frac{1-m-\epsilon}{(1-y)^3}\right],
\eea
then we have
\bea
\lambda_3+\lambda_4&=&\emph{tr}D-\mu=4\mu-2\delta\mu-\mu=-(2\delta-3)\mu,
\\
\lambda_3\lambda_4&=&-\mu(\lambda_3+\lambda_4)+{1\over2}(trD)^2-2\mu^2+\mu\sum_{i,j=1,i\ne j}^4 D_{ij}
  -{1\over2}\sum_{i=1}^4\left(\sum_{s=1,s\ne i}^4D_{is}\right)^2
  -\sum_{i,j=1,i<j}^4D_{ij}D_{ji}
\nonumber\\
&=&(2\delta-3)\mu^2+{1\over2}(4-2\delta)^2\mu^2-2\mu^2+2\delta\mu^2
-\frac{1}{2\alpha^6}\left[\left(\frac{\epsilon}{x^3}+\frac{\tau\epsilon}{y^3}+1-m-(1+\tau)\epsilon\right)^2\right.
\nonumber
\\
&&+\left(\frac{m}{x^3}+\frac{\tau\epsilon}{(y-x)^3}+\frac{1-m-(1+\tau)\epsilon}{(1-x)^3}\right)^2
+\left(\frac{m}{y^3}+\frac{\epsilon}{(y-x)^3}+\frac{1-m-(1+\tau)\epsilon}{(1-y)^3}\right)^2
\nonumber
\\
&&+\left(m+\frac{\epsilon}{(1-x)^3}+\frac{\tau\epsilon}{(1-y)^3}\right)^2
        +2\left(\frac{m\epsilon}{x^6}+\frac{m\tau\epsilon}{y^6}+m(1-m-(1+\tau)\epsilon)\right.
\nonumber
\\
&&\quad        \left.\left.+\frac{\tau\epsilon^2}{(y-x)^6}+\frac{\epsilon(1-m-(1+\tau)\epsilon)}{(1-x)^6}
               +\frac{\tau\epsilon(1-m-(1+\tau)\epsilon)}{(1-y)^6}\right)\right]
\nonumber\\
&=&\mu^2\Bigg[(2\delta^2-4\delta+3)-\frac{1}{2\mu^2\alpha^6}
\left[\left(\frac{\epsilon}{x^3}+\frac{\tau\epsilon}{y^3}+1-m-(1+\tau)\epsilon\right)^2\right.
\nonumber
\\
&&+\left(\frac{m}{x^3}+\frac{\tau\epsilon}{(y-x)^3}+\frac{1-m-(1+\tau)\epsilon}{(1-x)^3}\right)^2
+\left(\frac{m}{y^3}+\frac{\epsilon}{(y-x)^3}+\frac{1-m-(1+\tau)\epsilon}{(1-y)^3}\right)^2
\nonumber
\\
&&+\left(m+\frac{\epsilon}{(1-x)^3}+\frac{\tau\epsilon}{(1-y)^3}\right)^2
        +2\left(\frac{m\epsilon}{x^6}+\frac{m\tau\epsilon}{y^6}+m(1-m-(1+\tau)\epsilon)\right.
\nonumber
\\
&&\quad        \left.\left.+\frac{\tau\epsilon^2}{(y-x)^6}+\frac{\epsilon(1-m-(1+\tau)\epsilon)}{(1-x)^6}
               +\frac{\tau\epsilon(1-m-(1+\tau)\epsilon)}{(1-y)^6}\right)\right].
\eea
Moreover, we have
\bea
\Delta&=&(\lambda_3+\lambda_4)^2-4\lambda_3\lambda_4
\nonumber
\\
&=&\mu^2\left\{-4\delta^2+4\delta-3+\frac{2}{\mu^2\alpha^6}
\left[\left(\frac{\epsilon}{x^3}+\frac{\tau\epsilon}{y^3}+1-m-(1+\tau)\epsilon\right)^2\right.\right.
\nonumber
\\
&&+\left(\frac{m}{x^3}+\frac{\tau\epsilon}{(y-x)^3}+\frac{1-m-(1+\tau)\epsilon}{(1-x)^3}\right)^2
+\left(\frac{m}{y^3}+\frac{\epsilon}{(y-x)^3}+\frac{1-m-(1+\tau)\epsilon}{(1-y)^3}\right)^2
\nonumber
\\
&&+\left(m+\frac{\epsilon}{(1-x)^3}+\frac{\tau\epsilon}{(1-y)^3}\right)^2
        +2\left(\frac{m\epsilon}{x^6}+\frac{m\tau\epsilon}{y^6}+m(1-m-(1+\tau)\epsilon)\right.
\nonumber
\\
&&\quad        \left.\left.\left.+\frac{\tau\epsilon^2}{(y-x)^6}+\frac{\epsilon(1-m-(1+\tau)\epsilon)}{(1-x)^6}
               +\frac{\tau\epsilon(1-m-(1+\tau)\epsilon)}{(1-y)^6}\right)\right]\right\}
\eea
Letting $\tilde\Delta=\frac{\Delta}{4\mu^2}$,
and note that $\lambda_3\ge\lambda_4$ are real numbers, we have
\bea
\lambda_3&=&\left(-\delta+{3\over2}+\sqrt{\tilde\Delta}\right)\mu,
\\
\lambda_4&=&\left(-\delta+{3\over2}-\sqrt{\tilde\Delta}\right)\mu.
\eea
Therefore, we obtain
\be \label{bb12}
\bb_1=-\frac{\lambda_3}{\mu}=\delta-{3\over2}-\sqrt{\tilde\Delta},\quad
\bb_2=-\frac{\lambda_4}{\mu}=\delta-{3\over2}+\sqrt{\tilde\Delta}.\quad
\ee
Then using the numerical results by R. Mart\'{\i}nez, A. Sam\`{a} and C. Sim\'{o} in \cite{MSS1} and
\cite{MSS2}, we can obtain the stability pattern of our four body problem.

\medskip

{\bf Step 3.} {\it Computations on the limit case.}

We need to compute the mass parameter of the restricted three-body problem of given masses
$\tilde{m}_1=m$, $\tilde{m}_2=0$, and $\tilde{m}_3=1-m$.
By (A.3) of \cite{MSS2}, $\bb$ (they use $\bb_c$ there) is given by
\bea
\bb_c=-1+\frac{\alpha}{a^{\alpha+2}[1+(\rho+1)^2]}\left[(\rho+2)\frac{(\rho+1)m_1+m_2}{\rho^{\alpha+2}}
+(\rho+1)(m_2\rho+m_3(\rho+1))+\frac{m_3-m_1\rho}{(\rho+1)^{\alpha+2}}\right],\nonumber
\eea
where $\rho$ and $a$ is given by (A.2) of \cite{MSS2}.

Note that, when in our case $\alpha=1$, (A.2) of \cite{MSS2} is just the Euler's quintic equation,
then together with
$\rho=\frac{x_0}{1-x_0}$ of (\ref{quintic.polynomial}), we have
\bea
\beta&=&-1+\frac{(\rho+1)^3}{1+(\rho+1)^2}\left[(\rho+2)\frac{(\rho+1)m}{\rho^3}+(1-m)(\rho+1)^2+\frac{(1-m)-m\rho}{(\rho+1)^3}\right]
\nonumber\\
&=&-1+\frac{(\frac{x_0}{1-x_0}+1)^3}{1+(\frac{x_0}{1-x_0}+1)^2}\left[(\frac{x_0}{1-x_0}+2)\frac{(\frac{x_0}{1-x_0}+1)m}{(\frac{x_0}{1-x_0})^3}
+(1-m)(\frac{x_0}{1-x_0}+1)^2+\frac{(1-m)-m\frac{x_0}{1-x_0}}{(\frac{x_0}{1-x_0}+1)^3}\right]
\nonumber\\
&=&-1+\frac{1}{(1-x_0)[(1-x_0)^2+1]}\left[m\frac{(2-x_0)(1-x_0)}{x_0^3}+(1-m)\frac{1}{(1-x_0)^2}+(1-x_0)^2(1-x_0-m)\right]
\nonumber\\
&=&-1+\frac{1}{(1-x_0)[(1-x_0)^2+1]}\left[m\frac{(2-x_0)(1-x_0)}{x_0^3}-m\frac{(1-x_0)^2}{x_0^2}\right]
\nonumber\\
&&+\frac{1}{(1-x_0)[(1-x_0)^2+1]}\left[(1-m)\frac{1}{(1-x_0)^2}+(1-m)\right]
\nonumber\\
&&+\frac{1}{(1-x_0)[(1-x_0)^2+1]}\left[m\frac{(1-x_0)^2}{x_0^2}-(1-m)+(1-x_0)^2(1-x_0-m)\right]
\nonumber\\
&=&-1+\frac{m}{x_0^3}+\frac{1-m}{(1-x_0)^3}+\frac{1}{(1-x_0)[(1-x_0)^2+1]}\frac{-x_0^5+(3-m)x_0^4-(3-2m)x_0^3+mx_0^2-2mx_0+m}{x_0^2}.\label{mass.parameter}
\nonumber\\
&=&-1+\frac{m}{x_0^3}+\frac{1-m}{(1-x_0)^3},  \label{beta.expression}
\eea
where in the last equality, we used (\ref{x0Equation}).

Following pp.171 in \cite{Xia}, for $q=(q_x,q_y)^T\in\R^2$, we define
\be
V_2(q)=\frac{m}{|a_{1,0}-q|}+\frac{1-m}{|a_{4,0}-q|}+{1\over2}\alpha_0^{-3}|q|^2
\ee
where $\alpha_0^{-3}$ is an extra parameter because Z. Xia fixed $\lambda=1$ of (1) in \cite{Xia},
but here we have $\lambda=\alpha_0^{-3}$.
Then we have
\bea
\frac{\partial^2V_2}{\partial^2q_x}&=&-\frac{m}{|a_{1,0}-q|^3}-\frac{1-m}{|a_{4,0}-q|^3}+\frac{1}{\alpha_0^3}
      +3\left[\frac{m(-(1-m)\alpha_0-q_x)^2}{|a_{1,0}-q|^5}+\frac{(1-m)(m\alpha_0-q_x)^2}{|a_{4,0}-q|^5}\right],
\\
\frac{\partial^2V_2}{\partial q_x\partial q_y}
      &=&-3\left[\frac{m(-(1-m)\alpha_0-q_x)q_y}{|a_{1,0}-q|^5}+\frac{(1-m)(m\alpha_0-q_x)q_y}{|a_{4,0}-q|^5}\right]
\\
\frac{\partial^2V_2}{\partial^2q_x}&=&-\frac{m}{|a_{1,0}-q|^3}-\frac{1-m}{|a_{4,0}-q|^3}+\frac{1}{\alpha_0^3}
      +3\left[\frac{mq_y^2}{|a_{1,0}-q|^5}+\frac{(1-m)q_y^2}{|a_{4,0}-q|^5}\right].
\eea
Therefore
\bea
\frac{\partial^2V_2}{\partial^2q_x}\bigg|_{q=(x_0\alpha_0,0)^T}
&=&-\frac{m}{x_0^3\alpha_0^3}-\frac{1-m}{(1-x_0)^3\alpha_0^3}+\frac{1}{\alpha_0^3}
    +3\left[\frac{mx_0^2\alpha_0^2}{x_0^5\alpha_0^5}+\frac{(1-m)(1-x_0)^2\alpha_0^2}{(1-x_0)^5\alpha_0^5}\right]
\nonumber
\\
&=&{2\over{\alpha_0^3}}\left[\frac{m}{x_0^3\alpha_0^3}+\frac{1-m}{(1-x_0)^3\alpha_0^3}\right]+\frac{1}{\alpha_0^3}
\nonumber
\\
&=&(2\beta+3)\alpha_0^{-3},
\\
\frac{\partial^2V_2}{\partial q_x\partial q_y}\bigg|_{q=(x_0\alpha_0,0)^T}&=&0,
\\
\frac{\partial^2V_2}{\partial^2q_x}\bigg|_{q=(x_0\alpha_0,0)^T}
&=&-\frac{m}{x_0^3\alpha_0^3}-\frac{1-m}{(1-x_0)^3\alpha_0^3}+\frac{1}{\alpha_0^3}
\nonumber
\\
&=&-\beta\alpha_0^{-3},
\eea
and hence
\be
D^2V_2(q)|_{q=(x_0\alpha_0,0)^T}=\left(\matrix{(2\beta+3)\alpha_0^{-3}& 0\cr 0& -\beta\alpha_0^{-3}}\right).
\ee
By the Case (ii) in p.173 of \cite{Xia}, we have
\be
\lim_{\epsilon\rightarrow0}\frac{a_3-a_2}{(m_2+m_3)^{1\over3}}=\lim_{\epsilon\rightarrow0}r_2'=\pm[(2\beta+3)\alpha_0^{-3}]^{-{1\over3}},
\ee
and hence
\bea
\lim_{\epsilon\rightarrow0}\frac{m_2}{|a_2-a_3|^3}
&=&\frac{1}{1+\tau}\lim_{\epsilon\rightarrow0}\frac{m_2+m_3}{|a_2-a_3|^3}=\frac{(2\beta+3)\alpha_0^{-3}}{1+\tau}
=\frac{(2\beta+3)\mu_0}{1+\tau},
\label{m2_a23}
\\
\lim_{\epsilon\rightarrow0}\frac{m_3}{|a_2-a_3|^3}
&=&\frac{\tau}{1+\tau}\lim_{\epsilon\rightarrow0}\frac{m_2+m_3}{|a_2-a_3|^3}=\frac{\tau(2\beta+3)\alpha_0^{-3}}{1+\tau}
=\frac{\tau(2\beta+3)\mu_0}{1+\tau}.
\label{m3_a23}
\eea

Note that $m_2=\epsilon,m_3=\tau\epsilon$ and $\lim_{\epsilon\rightarrow0}|a_i-a_j|\ne0$ if $i<j,(i,j)\ne(2,3)$,
from (\ref{D_epsilon}), we have
\be\label{D_0}
D_0=\lim_{\epsilon\rightarrow0}D=
\left(
\matrix{
m\mu_0& 0& 0& (1-m)\mu_0\cr
\frac{m}{x_0^3}\mu_0&       [-\beta-\frac{\tau(2\beta+3)}{1+\tau}]\mu_0& \frac{\tau(2\beta+3)}{1+\tau}\mu_0&     \frac{1-m}{x_0^3}\mu_0&\cr
\frac{m}{x_0^3}\mu_0&    \frac{2\beta+3}{1+\tau}\mu_0& [-\beta-\frac{2\beta+3}{1+\tau}]\mu_0&     \frac{1-m}{x_0^3}\mu_0&\cr
m\mu_0& 0& 0& (1-m)\mu_0
}
\right),
\ee
where we have used (\ref{mu0}), (\ref{mass.parameter}), (\ref{m2_a23}) and (\ref{m3_a23}).
Then the characteristic polynomial of $D_0$ is given by
\bea
\det(D_0-\lambda I)&=&-\lambda(\mu_0-\lambda)
      \left|\matrix{
      [-\beta-\frac{\tau(2\beta+3)}{1+\tau}]\mu_0-\lambda& \frac{\tau(2\beta+3)}{1+\tau}\mu_0\cr
      \frac{2\beta+3}{1+\tau}\mu_0& [-\beta-\frac{2\beta+3}{1+\tau}]\mu_0-\lambda}
      \right|
\nonumber
\\
&=&\lambda(\lambda-\mu_0)(\lambda+\beta\mu_0)(\lambda+3(\beta+1)\mu_0).
\eea
Then all eigenvalues of $D_0$ are given by
\be\label{lambda1234}
\lambda_{1,0}=\mu_0,\quad\lambda_{2,0}=0,\quad\lambda_{3,0}=-\beta\mu_0,\quad\lambda_{4,0}=-3(\beta+1)\mu_0,
\ee
and hence by (\ref{bi}), we have
\bea
\bb_{1,0}&=&-\frac{\lambda_{3,0}}{\mu_0}=\beta,
\\
\bb_{2,0}&=&-\frac{\lambda_{4,0}}{\mu_0}=3(\beta+1).
\eea

From (\ref{lambda1234}), the four eigenvalues of $D_0$ are different,
then for $\epsilon>0$ small enough, we have
\be
\lim_{\epsilon\to0}\lambda_i=\lambda_{i,0},\quad 1\le i\le4.
\ee
Thus, we also have
\be
\lim_{\epsilon\to0}\bb_i=\bb_{i,0},\quad 1\le i\le2.
\ee

Therefore, the linear stability problem of the limiting case of our four-body problem when letting
$\ep\to 0$ is reduced to the linear stability problems of two restricted three-body problems, for
which one has mass parameter $\bb$, and the other has mass parameter $3(\bb+1)$. Then the numerical
results obtained by R. Mart\'{\i}nez, A. Sam\`{a} and C. Sim\'{o} in \cite{MSS1} and \cite{MSS2} can
be used to obtain the linear stability pattern of the limiting case of our four-body problem.
we will compute a concrete example at the end of this paper.

\medskip

\begin{example}
Computations on the actual case of the Earth-Moon-two space stations system.
\end{example}

We denote by ESSM system the short hand notation for the Earth-two space stations-Moon system.
From https://en.wikipedia.org/wiki/Earth and https://en.wikipedia.org/wiki/Moon, one can find
that the mass of Earth is $E=5.97237\times10^{24}$kg, the mass of the Moon is $M=7.342\times10^{22}$kg,
the distance between the Earth and the Moon is $d=384405$km, and the actual eccentricity of the orbit
of Moon is $e\approx0.0549$. This eccentricity is viewed as that of the orbits in the ESSM system.

By the normalization of the masses, we have
\be\label{related.mass.of.Earth}
m = \frac{E}{E+M} \approx 0.9879.
\ee
For two space stations in the line segment between the Earth and the Moon, as their masses tends to $0$
their limit position $x_0$ given by (\ref{x_0def}) is determined by (\ref{x0Equation}) and $m$. When
$m$ is given by (\ref{related.mass.of.Earth}), by a numerical computation, we have
\be
x_0 \approx 0.8493
\ee
By the distance between the Earth and the Moon, the distance between the limit position of the two
space stations and the Moon is $d_{SM}=d\times(1-x_0)\approx 57930$km.

Via (\ref{beta.expression}), the constant $\bb$ for the EEM of the $3$-body problem is given by
\be
\bb = -1+\frac{m}{x_0^3}+\frac{1-m}{(1-x_0)^3} \approx 4.1481.
\ee
Thus the linear stability property of the ESSM system is determined by the eccentricity
$e\approx0.0549$ of their orbits and the following two mass parameters:
\be
\bb_1 = \bb\approx4.1481,\qquad
\bb_2 = 3(\bb+1) \approx 15.4442.
\ee
On the other hand, by (1.5)-(1.8) of \cite{ZL}, we have
\be
\hat\bb_2 \approx 2.7122,\quad \hat\bb_{5\over2} \approx 4.9437,\quad
\hat\bb_{4} \approx 14.6764,\quad \hat\bb_{9\over2} \approx 18.9243,
\ee
where $\hat\bb_n$ and $\hat\bb_{n+{1\over 2}},n\in\N$
are the parameter values when the resonances of the linearized system appear.
Indeed, $\hat\bb_n$ is the $n$-th value such that $\ga_{\bb,0}(2\pi)$ has eigenvalue $1$,
and $\hat\bb_{n+{1\over2}}$ is the $n$-th value such that $\ga_{\bb,0}(2\pi)$ has eigenvalue $-1$.
Here $\ga_{\bb,0}(2\pi)$ is the end matrix at time $t=2\pi$ of the fundamental solution of the
linearized Hamiltonian system (\ref{2.41}) at the Euler solution EEM $q_{m,e}$ with $e=0$ of the
$3$-body problem. Hence in our case,
\be
\hat\bb_2<\bb_1<\hat\bb_{5\over2},\quad
\hat\bb_4<\bb_2<\hat\bb_{9\over2}.
\ee
Since the eccentricity $e \approx 0.0549$ is very small, numerical computations show that the
linear stability property is the same as that of $e=0$. Then by Theorem 1.5 of \cite{ZL},
the linear stability pattern of the ESSM system is
\be
R(\theta_1)\diamond D(2)\diamond R(\theta_2)\diamond D(2)
\ee
for some $\theta_1$ and $\theta_2\in (0,\pi)$. Here for $\th\in\R$ and $\lm\in \R\bs\{0,\pm 1\}$
we denote the elliptic and hyperbolic matrices by
$$ R(\th)=\left(\matrix{\cos\th & -\sin\th\cr
                        \sin\th  & \cos\th\cr}\right), \qquad
   D(\lm)=\left(\matrix{\lm & 0\cr
                         0  & \lm^{-1}\cr}\right), $$
respectively.

\setcounter{equation}{0}
\section{Appendix: A sketch of the proof of Lemma \ref{L1.1}.}

For reader's conveniences, following \cite{Moe1} of R. Moeckel (cf. also \cite{Pac}, \cite{Moe2}),
next we sketch the ideas of the proof of Lemma \ref{L1.1} due to C. Conley.

\medskip

{\bf A sketch of the proof of Lemma \ref{L1.1}.} Note first that both the matrices $D$ in (\ref{1.11})
and $\tilde{D}$ in (\ref{tilde.D}) possess the same eigenvalues by the definition of $\tilde{D}$. Because
$\tilde{D}$ is symmetric, all its eigenvalues are real, and then so does $D$, although it may not be
symmetric in general.

Note that $2B(a) = U''(a)$ is the Hessian of $U(q)$ at the collinear central configuration $q=a$,
and $U''(a)+U(a)\td{M}$ is the Hessian of $U|_{S}$ with $S$ being the hypersurface determined by
(\ref{1.7}). By the homogeneity of $U$, we obtain that $D$ has the first eigenvalue
$\lambda_1=\mu=U(a)$ with the eigenvector $v_1=(1,1,\ldots,1)^T$, i.e., $(Dv_1)_i=\mu$ holds for
$1\le i\le n$.

From the definition (\ref{eq.of.cc}) of $a$ as a central configuration, we obtain that $D$ has the
second eigenvalue $\lambda_2=0$ with the eigenvector $v_2=(a_{1x},a_{2x},\ldots,a_{nx})^T$. More
precisely for $1\le i\le n$ by (\ref{eq.of.cc}) we have
\bea (Dv_2)_i
&=& (\mu-\sum_{j=1,j\ne i}^n\frac{m_j}{|a_i-a_j|^3})a_{ix}+\sum_{j=1,j\ne i}^n\frac{m_ja_{jx}}{|a_i-a_j|^3} \nn\\
&=& \mu a_{ix}+\sum_{j=1,j\ne i}^n\frac{m_j(a_{jx}-a_{ix})}{|a_{jx}-a_{ix}|^3}  \nn\\
&=& 0.  \nn\eea

Note that by (\ref{1.6})-(\ref{1.7}), the vectors $v_1$ and $v_2$ form an $\td{M}$-orthonormal
sub-basis, i.e., they satisfy $v_1^T\tilde{M}v_1 = 1$, $v_1^T\tilde{M}v_2 =0$, and
$v_2^T\tilde{M}v_2 =1$. Denote all the other eigenvalues of $D$ by $\lm_3, \ldots, \lm_n$. Next
goal is to show that the other $(n-2)$ eigenvalues of $D$ are non-positive.

Following \cite{Moe1}, this is equivalent to showing that all the eigenvalues of $D$ are non-positive
when we restricted to the subspace spanned by vectors orthogonal to $\tilde{M}v_1$, and observing
that this is equivalent to showing that in the flow on the space of lines through the origin
determined by the following linear system on $u$,
\be \dot{u} = M^{-1}B(a)u, \label{surplus1}\ee
the line determined by $v_2$ is an attractor.

Let
$$ K = \left\{u=(u_1,u_2,\ldots,u_n)^T\;\left|\;
       \sum_{i=1}^nm_iu_i=0,\;u_1\le u_2\le\ldots\le u_n\right.\right\}.  $$
Then for any $u\in K$, we have $u\perp\tilde{M}v_1$. Moreover, we have $rv_2\in K$ for any $r\in\R$.
We will show that, around the line in $K$ which is carried strictly inside itself by the flow
defined by (\ref{surplus1}) except for the origin.

Note that the boundary $\partial K$ of $K$ consists of points where one or more equalities hold.
However, except for the origin, at least one strict inequality must hold, otherwise
$u=k(1,1,\ldots,1)^T\in K$ and hence $k=0$. Consider a boundary point with
$$  u_i = u_{i+1} = \cdots = u_j < u_{j+1},\quad 1\le i<j<n,  $$
or
$$  u_{i-1} = u_i = \cdots = u_j < u_{j+1},\quad 1<i<j\le n.  $$
The differential equation (\ref{surplus1}) becomes
$$ \dot{u}_i=\sum_{k\ne i}{m_k\over r_{ik}^3}(u_k-u_i),\quad
   \dot{u}_j=\sum_{k\ne j}{m_k\over r_{jk}^3}(u_k-u_j),  $$
where $r_{mk}=|a_{mx}-a_{kx}|$ for $m=i, j$. Since $u_i=u_j$ we get
$$ \dot{u}_j-\dot{u}_i
    = \sum_{k\ne i,j}m_k(u_k-u_j)\left[{1\over r_{jk}^3}-{1\over r_{ik}^3}\right]. $$
Every term in this sum is non-negative, since

(i) if $k<i$, $\;(u_k-u_i)\le0$ and ${1\over r_{jk}^3}-{1\over r_{ik}^3}<0$;

(ii) if $i\le k\le j$, $\;u_k-u_i=0$;

(iii) if $k>j$, $\;(u_k-u_i)\ge 0$ and ${1\over r_{jk}^3}-{1\over r_{ik}^3}>0$.

Moreover, at least one term is strictly positive since not all of $u_i$ with $1\le i\le n$
are equal. Thus $\dot{u}_j-\dot{u}_i>0$ and the boundary point moves into the interiors of
the cone $K$ as required.

Now we consider the central configurations in $\R^3$. Let
\bea
S &=& \left\{q=(q_1,q_2,\ldots,q_n)^T,q_i\in\R^3\;\Bigg|\;
      \sum_{i=1}^nm_iq_i=0,\sum_{i=1}^nm_iq_i^2=1,q_i\ne q_j\;{\rm if}\;i\ne j\right\}, \nn\\
C &=& \left\{q=(q_1,q_2,\ldots,q_n)^T\in S\;\Bigg|\;
      q_i\in\R\times\{0\}\times\{0\},\quad\forall 1\le i\le n\right\}, \nn\\
E &=& \left\{q=(q_1,q_2,\ldots,q_n)^T\in S\;\Bigg|\;
      q_i\in\{0\}\times\R^2,\quad\forall 1\le i\le n\right\}, \nn\\
\tilde{C} &=& \left\{q=(q_1,q_2,\ldots,q_n)^T\in S\;|\; q\; {\rm is\;collinear\;
      along\; some\; line}\right\}.  \nn\eea
Then $C\subset\tilde{C}$ holds and $\tilde{C}$ is the orbit of $C$ under $\SO(3)$.

Now on $S$, the central configuration equation is
$$ F(q)=\tilde{M}^{-1}U'(q)+U(q)q=0, $$
where $U'(q)$ denotes the gradient of $U$ with respect to $q=(q_1,\ldots, q_n)$. Then when
we consider the gradient flow of the system
\be \dot{q}= F(q), \lb{gradient.flow}\ee
a central configuration is a fixed point of this flow. Note that $C$, $\tilde{C}$ and $E$
are invariant sub-manifolds under the gradient flow of (\ref{gradient.flow}). For the central
configuration $q_0=(q_{1,0},q_{2,0},\ldots,q_{n,0})$ with $q_{i,0}=(a_{ix},0,0)^T$, we have
\bea
F'(q_0)|_C &=& -2\tilde{M}^{-1}B+\mu I_n,  \\
F'(q_0)|_E &=& \diag(\tilde{M}^{-1}B, \tilde{M}^{-1}B)+\mu I_{2n}.  \label{DF.E}\eea

Note that in the first Corollary on p.507 of \cite{Moe1}, R. Moeckel proved that any orbits near
$\tilde{C}$ are attracted to $\tilde{C}$ by the gradient flow of (\ref{gradient.flow}). Therefore
it yields that $F'(q_0)|_E$ in (\ref{DF.E}) is non-negative definite as required. In fact, using
notations in \cite{Moe1}, an explicit neighborhood
$\mathcal{U}=\{q\in S|\;\Theta(q)\le{\pi\over4}\}$ of $\tilde{C}$ in $S$ can be defined such that
the orbits of the gradient flow of (\ref{gradient.flow}) in $\mathcal{U}$ get more and more
collinear.

Here following \cite{Moe1} the function $\Theta(q,L)$ measures the approximate collinearity of
a configuration $q\in S$ and a line $L$ in $\R^3$ is defined by
$$  \Theta(q,L)=\max_{i\ne j}\angle(L,q_i-q_j),  $$
where $\angle(L,q_i-q_j)$ denotes the acute angle between $L$ and $q_i-q_j$. $\Theta(q,L)$
vanishes if and only if $q$ is collinear along a line parallel to $L$. Then let
$$ \Theta(q)=\min_{L}\Theta(q,L), $$
which vanishes if and only if $q$ is collinear.

Note that in $\mathcal{U}$, $\Theta(q)$ is strictly decreasing along orbits $q=q(t)$ of the
gradient flow of (\ref{gradient.flow}), and it suffices to prove
\be \Theta(q(t))<\Theta(q(0)), \qquad \forall\;t>0. \lb{A1}\ee
Now we refer readers to pp.504-505 of \cite{Moe1} on the details of the proof of (\ref{A1}). \hb

\medskip

\noindent {\bf Acknowledgements.} The authors would like to thank sincerely the anonymous editor for
informing us and helps on finding the paper of J. Liouville, and valuable comments. They thank sincerely
also the anonymous referees on their careful reading and helpful comments on the manuscript of
this paper.


\begin{thebibliography}{}


\bibitem{And} M. Andoyer, Sur les solutiones periodiques voisines des position
d'equilibre relatif dans ie probleme des n corps. Bull. Astron. 23, (1906) 129-146.

\bibitem{Dan} J. Danby, The stability of the triangular Lagrangian point in the general
problem of three bodies. {\it Astron. J.} 69. (1964) 294-296.

\bibitem{Euler} L. Euler, De motu restilineo trium corporum se mutus
attrahentium. {\it Novi Comm. Acad. Sci. Imp. Petrop.}  11. (1767) 144-151.

\bibitem{Ga} M. Gascheau,  Examen d'une classe d'\'{e}quations diff\'{e}rentielles
et application \`{a} un cas particulier du probl\`{e}me des trois corps. {\it Comptes
Rend. Acad. Sciences.} 16. (1843) 393-394.

\bibitem{HLS} X. Hu, Y. Long, S. Sun, Linear stability of elliptic Euler
solutions of the classical planar three-body problem via index theory.
{\it Arch. Ration. Mech. Anal.} 213. (2014) 993-1045.

\bibitem{HO} X. Hu, Y. Ou, Collision index and stability of elliptic relative equilibria
in planar $n$-body problem. {\it http://arxiv.org/pdf/1509.02605.} (2015). {\it Comm.
Math. Phys.} to appear.

\bibitem{HS} X. Hu, S. Sun, Morse index and stability of elliptic Lagrangian
solutions in the planar three-body problem. {\it Advances in Math.} 223. (2010) 98-119.

\bibitem{IM} R. Iturriaga, E. Maderna, Generic uniqueness of the minimal
Moulton central configuration. {\it http://arxiv.org/abs/1406.6887v3.} (2015).

\bibitem{Lag} J. Lagrange,  Essai sur le probl\`{e}me des
trois corps. Chapitre II. {\OE}uvres Tome {6}, Gauthier-Villars,
Paris. (1772) 272-292.

\bibitem{Lio} J. Liouville, Sur un cas particulier du probl\`{e}me des trois corps.
{\it J. Math. Pures Appl.} 7. (1842) 110-113.

\bibitem{Lon} Y. Long, Lectures on Celestial Mechanics and Variational Methods.
{\it Preprint.} 2012.

\bibitem{MSS} R. Mart\'{\i}nez, A. Sam\`{a}, C. Sim\'{o},
Stability of homograpgic solutions of the planar three-body problem
with homogeneous potentials. in International conference on
Differential equations. Hasselt, 2003, eds, Dumortier, Broer, Mawhin,
Vanderbauwhede and Lunel, World Scientific, (2004) 1005-1010.

\bibitem{MSS1} R. Mart\'{\i}nez, A. Sam\`{a}, C. Sim\'{o},
Stability diagram for 4D linear periodic systems with applications
to homographic solutions. {\it J. Diff. Equa.} 226. (2006) 619-651.

\bibitem{MSS2} R. Mart\'{\i}nez, A. Sam\`{a}, C. Sim\'{o}, Analysis of
the stability of a family of singular-limit linear periodic systems in
$\R^4$. Applications. {\it J. Diff. Equa.} 226. (2006) 652-686.

\bibitem{Mey} M. Meyer, Solutiones voisines des solutiones de lagrange dans Ie
probleme des n corps. Ann. Obs. Bordeaux, 17 (1933) 77-252.

\bibitem{MS} K. Meyer, D. Schmidt, Elliptic relative equilibria in
the N-body problem. {\it J. Diff. Equa.} 214. (2005) 256-298.

\bibitem{Moe1} R. Moekel, On central configurations. {\it Math. Z.} 205 (1990) 499-517.

\bibitem{Moe2} R. Moekel, Celestial Mechanics (especially central configurations).
{\it http://www.math.umn.edu/~rmoeckel/notes/CMNotes.pdf}. 1994.

\bibitem{Moe3} R. Moekel, Linear stability analysis of some symmetrical classes
of relative equilibria. H. S. Dumas et al. (eds.), Hamiltonian Dynamical Systems,
Springer, New York. (1995) 291-317.

\bibitem{Mou} F. Moulton, The straight line solutions of the $n$-body problem.
{\it Ann. of Math}. II Ser. 12 (1910) 1-17.

\bibitem{Pac} F. Pacella, Central configurations and the equivariant Morse theory.
{\it Arch. Ration. Mech. Anal.} 97 (1987) 59-74.

\bibitem{R1} G. Roberts, Linear stability of the elliptic Lagrangian
triangle solutions in the three-body problem. {\it J. Diff. Equa.}
182. (2002) 191-218.

\bibitem{R2} E. Routh, On Laplace's three particles with
a supplement on the stability or their motion. {\it Proc. London Math.
Soc.} 6. (1875) 86-97.

\bibitem{Xia} Z. Xia, Central Configurations with Many Small Masses.
{\it J. Diff. Equa.} 91. (1991) 168-179.

\bibitem{Win1} A. Wintner, The Analytical Foundations of Celestial Mechanics.
Princeton Univ. Press, Princeton, NJ. 1941. Second print, Princeton Math. Series
5, 215. 1947.

\bibitem{ZL} Q. Zhou, Y. Long, Maslov-type indices and linear stability
of elliptic Euler solutions of the three-body problem.
{\it http://arxiv.org/pdf/1510.06822v1.} (2015). {\it Submited.}


\end{thebibliography}
\end{document}